\documentclass[a4paper,12pt]{article}
\usepackage{hyperref}
\usepackage{amssymb}
\usepackage{amsmath}
\usepackage[latin1]{inputenc} 
\newtheorem{defi}{Definition}[section]
\newtheorem{prop}[defi]{Proposition}
\newtheorem{thm}[defi]{Theorem}
\newtheorem{lem}[defi]{Lemma}
\def\R{{\mathbb R}}
\def\C{{\mathbb C}}
\def\Z{{\mathbb Z}}
\def\N{{\mathbb N}}

\def\F{{\mathbb F}}
\def\OO{{\mathcal O}}
\def\pro{\mathrm{p}}
\def\CC{{\mathcal C}}
\def\E{\mathcal{E}}

\def\eps{{\varepsilon}}
\def\fixme{{$ $}}
\def\pr{{\mathrm{p}}}
\def\lie{{\mathrm{Lie}}}

\newcommand {\tx}[1] {\textrm{#1}}
\newcommand {\resr}[1] {\OO/\pi^{#1}\OO}
\newcommand {\quotr}[1] {\OO/\pi^{#1}\OO}
\newcommand {\esp}[1]{\underset{ #1}{\mathbb E}}
\newcommand{\qed}{\hfill $\Box$}
\newcommand{\paire}[2]{\langle{#1},{#2}\rangle}
\newcommand{\flf}[1]{\lfloor{#1}\rfloor}
\newcommand{\clf}[1]{\lceil{#1}\rceil}

\begin{document}

\title{Strong Banach Property (T) for Simple Algebraic Groups of Higher
Rank}
\author{Benben Liao\footnote{Institut de mathématiques de Jussieu (Université
Paris Diderot- Paris 7).}} 
\date{\today}
\maketitle
\abstract{In \cite{laf-duke,laf-jta}, Vincent Lafforgue proved strong Banach property (T)
for $SL_3$ over a non archimedean local field $F.$ In this paper, we extend his
results to $Sp_4$ and therefore to any connected almost $F$-simple algebraic
group with $F$-split rank $\geq 2.$ As applications, the family of expanders constructed by finite quotients of
a lattice in such a group does not admit a uniform embedding
in any Banach space of type $>1,$ and any affine isometric action of such a group, or of any cocompact lattice in it,
in a Banach space of type $>1$ has a fixed point.}
\tableofcontents

\section{Introduction} 

In \cite{laf-duke,laf-jta}, Vincent Lafforgue proved strong Banach property (T)
for $SL_3$ over a non archimedean local field $F.$ In this paper, we extend his
results to $Sp_4$ and therefore to any connected almost $F$-simple algebraic
group with $F$-split rank $\geq 2.$ As the first application, the family of expanders constructed
by finite quotients of a lattice in such a group does not admit a uniform embedding
in any Banach space of type $>1.$ As the second application, we prove that 
any affine isometric action of such a group, or of any cocompact lattice in it, in a Banach space of type $>1$ has a fixed point.
In \cite{bader}, it is conjectured that any isometric affine action of a higher rank simple algebraic group over a local field and of its lattice in a uniformly convex space has a fixed point. As a consequence of the second application, we confirm this conjecture for any non archimedean local field and the corresponding cocompact lattices.

To announce the precise statements, we begin by recalling some definitions and
notations from \cite{laf-jta}. 
\begin{defi}\label{type}
A class of Banach spaces $\E$ is of type $>1$ if one of the following two
equivalent conditions holds.
\begin{itemize}
  \item i) There exist $n\in\N$ and $\eps>0$ such that for any Banach space
  $E\in\E,$ $E$ does not contain $\ell_1^n$ $(1+\eps)$-isometrically;
  \item ii) There exist $p>1$ (called the type) and $T \in \R_{+}$ such that
  for any $E\in\E$, $n \in \N^{*}$ and $x_{1},...,x_{n} \in E$, we have
  $$\Big(\esp{\eps_{i}=\pm 1} \|\sum
  _{i=1}^{n}\eps_{i}x_{i}\|_{E}^{2}\Big)^{\frac{1}{2}}\leq T\Big(
  \sum_{i=1}^{n} \|x_{i}\|^{p}_{E}\Big)^{\frac{1}{p}}.$$
\end{itemize}
\end{defi} 
{\bf Remark 1.} We say that a class of Banach spaces $\E$ is given by
a super-property, if any Banach space $F$ finitely representable in $\E$ (i.e. for any finite dimensional
subspace $V\subset F$ and $\eps>0$ there exists $E\in\E$ which contains $V$ $(1+\eps)$-isometrically)
is an element of $\E.$ It is clear that a class of type $>1$ is given by a
super-property.
\\{\bf Remark 2.} If $\E$ is a class of Banach spaces given by a super-property
and not a class of type $>1,$ then $\E$ contains $L_1(\mu),$
where $\mu$ is any $\sigma$-finite measure. In fact, by the classification
of $\sigma$-finite measures it suffices to show that $\ell_1$ and
$L_1(\{0,1\}^\infty)$ are elements of $\E.$ $L_1(\{0,1\}^\infty)$ is finitely
representable in $\ell_1.$ By condition i) in the definition, $\ell_1$ is
finitely representable in the class $\E.$ Since $\E$ is given by a
super-property, we conclude that $L_1(\{0,1\}^\infty)$ and $\ell_1$ belong to
$\E.$

Let $\E$ be a class of Banach spaces stable under
complex conjugation and duality. Let $G$ be a locally compact topological group.
Let $\ell$ be a continuous length function of $G.$
Denote $\E_{G,\ell}$ the set of isomorphism classes of strongly continuous representations $(E,\pi)$
of $G$ such that $E\in\E$ and $$\|\pi(g)\|_{\mathcal L(E)}\leq e^{\ell(g)}$$ for
any $g\in G.$ Denote $\CC^\E_{\ell}(G)$ the completion of compactly supported
functions $C_c(G)$ on $G$ with respect to the norm
$$\|f\|_{\CC^\E_{\ell}(G)}=\tx{sup}_{(E,\pi)\in\E_{G,\ell}}\|\int
f(g)\pi(g)dg\|_{\mathcal L(E)}.$$

\begin{defi}\label{strong-ban-T} 
We say that a locally compact group $G$ has strong Banach property (T) if for
any class of Banach spaces $\E$ of type $> 1,$ stable under complex conjugation
and duality, and any continuous length function $\ell$ over $G$, there exists $s_0>0$ such
that the following holds. For any $C>0$ and $s_0\geq s\geq 0$, there exists a
real self-adjoint idempotent element $\pro$ in $\CC^\E_{C+s\ell}(G)$, such that
for any representation $(E,\pi)\in\E_{G,C+s\ell},$ the image of $\pi(\pro)$
consists of all $G$-invariant vectors in $E,$ i.e.
$$\pi(\pr)E=E^{\pi(G)}.$$
\end{defi}

{\bf Remark. }In this definition, the condition of type $>1$ cannot be replaced
by a weaker condition given by a super-property because otherwise it would be
satisfied only for compact groups.
Indeed when $G$ is non
compact, suppose that $\E$ is a class of Banach spaces (stable under complex
conjugation and duality) given by a super-property, and that there exists a real
self-ajoint idempotent $\pr\in C_{0}^\E(G)$ such that for any $(E,\pi)\in\E_{G,0}$ we have $\pi(\pr)E=E^{\pi(G)}$, we show that 
$\E$ is a class of Banach spaces of
type $>1.$ If not, by remark 2 below definition \ref{type},
$\E$ must contain $L^1(G).$ Note that for any
$(E_1,\pi_1),(E_2,\pi_2)\in\E_{G,0}$, any surjective morphism $E_1\to E_2$ in
the category $\E_{G,0}$ induces a surjective morphism from $E_1^G=\pi_1(\pr)E_1$
to $E_2^G=\pi_2(\pr)E_2.$ Now consider the morphism from $L^1(G)$ (with the left regular representation of
$G$) to $\C$ (with the trivial action of $G$) by integration on $G.$ Since $G$
is non compact, there is no non zero $G$-invariant integrable function on $G$,
therefore $L^1(G)^G=\{0\}$. However, $\C^G=\C$, and
this is a contradiction to that $L^1(G)^G\to \C^G$ must be a surjective
morphism. Therefore, $\E$ must be a class of type $>1$ (see the second remark
below definition 0.2 in \cite{laf-jta}).

Let $F$ be a non archimedean local field. The purpose
of this paper is to prove the following theorem.

\begin{thm}\label{main-thm}
Any connected almost $F$-simple algebraic group with $F$-split rank
$\geq 2$ has strong Banach property (T).
\end{thm}

{\bf Remark.} This result cannot be extended to any almost $F$-simple algebraic
group with $F$-split rank $= 1$ because they do not even have Kazhdan's
property (T). 

The following definition corresponds to the special case of isometric actions.
\begin{defi}\label{ban-T}
We say that a locally compact group $G$ has Banach property (T) if for any
class of Banach spaces $\E$ of type $> 1$ stable under complex conjugation
and duality, there exists a
real self-adjoint idempotent element $\pro$ in $\CC^\E_{0}(G)$, such that
for any representation $(E,\pi)\in\E_{G,0},$ the image of $\pi(\pro)$
consists of all $G$-invariant vectors in $E.$
\end{defi}
{\bf Remark. } If a locally compact group $G$ has (strong) Banach property (T)
with $\pr\in\CC^\E_{C+s\ell}(G)$ being the corresponding idempotent, there
always exist $\pr_n\in C_c(G)$ of integral $1,$ such that $\pr_n$ converges to
$\pr$ in $\CC^\E_{C+s\ell}(G).$ 
In fact, let $\check \pro_n\in C_c(G)$ be any sequence such that $\check
\pro_n\to \pro.$ Let $s_n=\int_G\check\pro_n(g) dg.$ Then 
\begin{gather*}\|\pr-s_n\pr\|_{\CC^\E_{C+s\ell}(G)}=\|\pr^2-\check\pr_n\pr\|_{\CC^\E_{C+s\ell}(G)}\\
\leq\|\pr-\check\pr_n\|_{\CC^\E_{C+s\ell}(G)}\|\pr\|_{\CC^\E_{C+s\ell}(G)},\end{gather*}
and hence $|1-s_n|\leq\|\pr-\check\pr_n\|_{\CC^\E_{C+s\ell}(G)}\to 1$ when
$n\to\infty.$ Therefore, $s_n\neq 0$ for big enough $n$ and
$\pro_n=\check\pro_n/s_n$ has integral $1$ and tends to $\pro.$

With the remark above and the same argument as in theorem 5.4 in \cite{laf-jta}, we
obtain the following theorem \ref{expanders} on application to expanders. 

We say that
 a family of graphes $\{(X_{i},d_i)\}_{i\geq 1}$ is embedded uniformly in a Banach space $E$, if there exist
 a function $\rho:\N\to\R_{+}$ that tends to infinity at infinity and
 $1$-Lipschitz maps $f_{i}:X_{i}\to E$ such that
 $$\|f_{i}(x)-f_{i}(y)\|_{E}\geq\rho(d_{i}(x,y))$$ for any $i\in\N$ and $x,y\in X_{i}.$ 
 
Let
$\Gamma$ be a discrete group with Banach property (T). Let $(\Gamma_{i})_{i\in\N}$ be a family of subgroups of $\Gamma$
such that $|\Gamma/\Gamma_{i}|$ tends to infinity. Let $S$ a finite symmetric
 system of generators of $\Gamma$ which contains $1.$
 For any $i\geq 0,$ $X_{i}=\Gamma/\Gamma_{i}$ is endowed with a graph structure
 associated to $S$ and we denote by $d_{i}$ the associated metric. As $\Gamma$
 has the usual property (T), $X_{i}$ forms a family of expanders. 
 
\begin{thm}\label{expanders}
Let $\Gamma$ be any discrete group
with Banach property (T). Then the family of expanders $(X_{i},d_{i})$
constructed above does not admit a uniform embedding in any Banach space of type
$>1$.
\end{thm}
Since strong Banach property clearly implies Banach property (T), and Banach
property (T) is inherited by lattices (proposition 5.3 in \cite{laf-jta}),
when $\Gamma$ is a lattice of a connected almost $F$-simple algebraic groups
of $F$-split rank $\geq 2,$ we
see that the family of expanders constructed above does not admit a
uniform embedding in any Banach space of type $>1.$

We recall that it is still unknown whether or not such a family of expanders
(or in fact any family of expanders) admits a uniform embedding in a Banach of
finite cotype (see \cite{laf-jta}, \cite{pis} and \cite{mendel-naor}).

We turn to application to fixed-point property. As a consequence of proposition
5.6 in \cite{laf-jta}, we immediately obtain
the following proposition, confirming conjecure 1.6 in \cite{bader} for any simple algebraic group of higher rank over a non archimedean local field and its cocompact lattice.

\begin{prop}
Let $G$ be a connected almost
$F$-simple algebraic group with $F$-split rank $\geq 2,$ or a cocompact lattice of such a group. Then any
affine isometric action of $G$ on a Banach space of type $>1$
has a fixed point.
\end{prop}
{\bf Remark 1.} This result cannot be strengthened to affine isometric actions for a larger class of Banach spaces defined by a super-property. If so, first of all by remark 2 below definition \ref{type} this class must contain all $L^1$ spaces and their closed subspaces. Denote $d\mu$
the Haar measure on $G,$ and $L^1_{i}(G)$ the space of functions $f\in L^1(G)$ such that $\int_G
f(g)d\mu(g)=i, i=0,1.$ Then $L^1_1(G)$ is an affine Banach space with $L^1_0(G)$
as the underlying Banach space. Let $G$ act on $L_1^1(G)$ by left translation. 
It is an affine isometric action of $G$ without fixed point, since $G$ is not
compact.
\\{\bf Remark 2.} As pointed out by Mikael de la Salle and the editor, let us mention that it is shown in \cite{monod11} that fixed point property for all $L^1$ spaces is a characterization of Kazhdan's property (T) for locally compact topological groups.

This paper will be part of my PhD thesis in Université Paris Diderot-
Paris 7. I would like to thank my thesis adviser Vincent Lafforgue for his
encouragement and guidance, and very helpful discussions about this paper. I
also thank Yanqi Qiu for the discussion of type of a Banach space.

Here is how the paper is organized. In section 2, we review the
theorem of strong Banach property (T) for $SL_3$ in \cite{laf-jta} and announce
the corresponding theorem
\ref{thm-sp4} for $Sp_4.$ In section 3, we prove theorem \ref{thm-sp4}
when $\tx{char}(F)\neq 2$ by constructing matrices for $Sp_4$ and
adapting the arguments in \cite{laf-jta}. In section 4,
we prove theorem \ref{thm-sp4} when $\tx{char}(F)=2$ by constructing new matrices for
the local estimate of the move $(0,2)$ and establishing the existence of two
limits in the spherical proposition.
In section 5,
we adapt a well known argument \cite{dk,vas,wang} and extend the results of
$SL_3$ and $Sp_4$ to any almost $F$-simple algebraic groups with $F$-split rank
$\geq 2.$ 

\section{Strong Banach property (T) for $Sp_4(F)$}\label{section-sp4}
Let $\E$ be any class of Banach spaces of type $>1$, stable under complex
conjugation and duality. Let $F$ be a non archimedean local field, $\OO$ the ring of
integers of $F$, $\pi$ one of
its uniformizer, $\F$ the residue field, and $q$ the cardinality of $\F,$ i.e. $q=\frac{1}{|\pi|}$.
The following proposition from \cite{laf-jta} (corollary 2.3) introduces
parameters $\alpha>0$ and $h\in\N^*$ for the class $\E.$
\begin{prop}\label{fourier-type}
There exist $\alpha
>0$ and $h\in\N^*$ such that for any $E\in\E$ we have $$\|T_{\resr{h}}\otimes
1_E\|\leq
e^{-\alpha},$$ where $T_{\quotr{h}}\otimes 1_E\in \mathcal
L\big(\ell^2(\quotr{h},E),\ell^2(\widehat {\quotr{h}},E)\big)$ is defined by
$$\big(T_{\quotr{h}}\otimes 1_E\big)(f)(\chi)=\esp{a\in\quotr{h}}\chi(a)f(a),$$
for any $\chi\in\widehat{\quotr{h}}$ and $f\in\ell^2(\quotr{h},E).$
\end{prop}

It is proved in \cite{laf-jta} that $SL_3(F)$ has strong Banach property
(T). 
\begin{thm}(Theorem 4.1 of \cite{laf-jta})\label{thm-sl3} Let $G=SL_3(F),$ and
$\ell$ be the length function on $G$ defined by
$$\ell\Big(k(\pi^{\frac{i+2j}{3}}\begin{pmatrix}\pi^{-i-j}&&\\&\pi^{-j}&\\&&1\end{pmatrix})k'\Big)=i+j,$$
for any $k,k'\in SL_3(\OO)$ and $i,j\geq 0$ with $i- j\in3\Z.$ Let $\beta
\in [0,\frac{\alpha}{3h})$. There exist $t,C'>0$ such that for any $C \in \R_{+}$, there exists a real and self-adjoint idempotent element $\pro\in \CC^{\E}_{C+\beta\ell}(G)$ such that \begin{itemize} 
	\item (i) for any representation $(E,\pi) \in \E_{G,C+\beta \ell }$, the image
	of $\pi(\pro)$ is the subspace of E consisting of all $G$-invariant vectors, 
	\item (ii) there exists a sequence $\pro_{n}\in C_{c}(G)$, such that
	$\int_{G}|\pro_{n}(g)|dg\leq 1$, $\pro_{n}$ has support in $\{g\in
	G,\ell(g)\leq n\}$, and $$\|\pro-\pro_{n}\|_{\CC^{\E}_{C+\beta\ell}(G)}\leq C'
	e^{2C-tn}.$$
\end{itemize}
\end{thm}

 Now we turn to $Sp_4.$ Let
$G=Sp_4(F)$, which is the group of $4\times 4$ matrices $g$ over $F$ such that
$^tgJg=J$ where $J$ is the skew-symmetric matrix, 
$$J=\begin{pmatrix}0&0&0&1\\0&0&1&0\\0&-1&0&0\\-1&0&0&0\end{pmatrix}.$$ Let
$K=Sp_4(\OO)$ (i.e. the subgroup in $Sp_4(F)$ whose matrix elements are in
$\OO$).
For any $(i,j)\in\Z^2$ let $$D(i,j)=\begin{pmatrix}\pi^{-i} & & & \\
 & \pi^{-j} & & \\
 & & \pi^j & \\
 & & & \pi^i
\end{pmatrix}.$$  
By $\|g\|$ we denote the norm of the operator $g\in \text{End}(F^4)$ w.r.t. the
standard norm on $F^4,$ i.e. $\|g\|=\max_{1\leq \alpha,\beta\leq
4}|g_{\alpha\beta}|.$ Similarly, denote $\|\Lambda^2g\|$ the biggest norm of all
$2\times 2$ minors of $g\in G,$ which is the norm of $\Lambda^2
g\in\text{End}(\Lambda^2F^4)$ w.r.t. the standard norm on $\Lambda^2F^4.$
 Let $\Lambda=\{(i,j)\in\N^2, i\geq j\}.$ Any element in $G$ has the form
 $kD(i,j)k'$ for some $(i, j)\in\Lambda$ and $k,k'\in K.$ For
such a $g=kD(i,j)k'\in G,$ we have $\|g\|=q^i$ and $\|\Lambda^2g\|=q^{i+j},$ and
this gives a bijection from $K\backslash G/K$ to $\Lambda$ by $g\mapsto(i,j),$ which is the
inverse of $(i,j)\mapsto KD(i,j)K.$ Let $\ell$ be
the length function of $G$ defined by $\ell(kD(i,j)k')=i+j,$ for any $k,k'\in K$ and
$(i,j)\in\Lambda$.

We will prove the following theorem with the argument used in \cite{laf-jta} for
the proof of theorem \ref{thm-sl3} (note that the statement is the same except
for the range of $\beta$).

\begin{thm}\label{thm-sp4}
 Let $\alpha$ and $h$ be as in proposition \ref{fourier-type}, and $\beta \in
 [0,\frac{\alpha}{8h})$.
 There exist $t,C'>0$ such that for any $C \in \R_{+}$, there exists a real and self-adjoint idempotent element $\pro\in
\CC^{\E}_{C+\beta\ell}(G)$ such that \begin{itemize} 
	\item (i) for any representation $(E,\pi) \in \E_{G,C+\beta \ell }$, the image
	of $\pi(\pro)$ is the subspace of E consisting of all $G$-invariant vectors, 
	\item (ii) there exists a sequence $\pro_{n}\in C_{c}(G)$, such that
	$\int_{G}|\pro_{n}(g)|dg\leq 1$, $\pro_{n}$ has support in $\{g\in
	G,\ell(g)\leq n\}$, and $$\|\pro-\pro_{n}\|_{\CC^{\E}_{C+\beta\ell}(G)}\leq C'
	e^{2C-tn}.$$ \end{itemize}
\end{thm}

\section{Proof of theorem \ref{thm-sp4} when $\tx{char}(F)\neq 2$}
This section is dedicated to the proof of theorem \ref{thm-sp4} when the
characteristic of $F$ is different from $2.$ We will first reduce the theorem
to two propositions on matrix
coefficients, and then prove them by
a zig-zag argument in the Weyl
chamber with two local estimates
of the matrix coefficients.

Most of the claims in this
section are only true when $\tx{char}(F)\neq 2,$ but some are still valid in
characteristic $2$ and will be used in the next section for the proof in
characteristic $2.$ 

When $\tx{char}F\neq 2,$ we denote $v_0$ the valuation of
$2\in\OO.$ For any $a\in\R,$ denote $\flf{a}$
(resp. $\clf{a}$) the biggest (resp. smallest) integer $\leq a$ (resp. $\geq
a$).

 
 Let
$(E,\pi)$ be any continuous representation of $G$ of a Banach space $E,$
$(V,\tau)$ any irreducible unitary representation of $K.$ For fixed $\xi\in
E$ and $\eta \in V\otimes E^*,$ we denote
$c(g)=\paire{\eta}{\pi(g)\xi}\in V$ for any $g\in G.$ By abuse of notation we
write $$c(i,j)=\paire{\eta}{\pi(D(i,j))\xi}.$$

The following is the proposition on spherical matrix coefficients, which will be
used to construct the idempotent element $\pr$ in theorem \ref{thm-sp4}.
\begin{prop}\label{spher-prop}
Suppose that $\tx{char}(F)\neq 2.$ Let $\alpha$ be as in proposition
\ref{fourier-type}, $\beta\in[0, \frac{\alpha}{4h})$.
There exists $C'>0$, such that the following holds. Let $C\in\R_+^*$, $(E,\pi)$
any element in $\mathcal E_{G,C+\beta\ell}$, and $\xi\in E$, $\eta\in E^*$ any
$K$-invariant vectors of norm $1$. 
There exists $c_\infty\in\C$, such that for
any $i\geq j\geq 0,$$$|c(i,j)-c_\infty|\leq
C'e^{2C-(\frac{\alpha}{2h}-2\beta)i}.$$
\end{prop}

Next we turn to the proposition on non spherical matrix coeffients.

\begin{prop}\label{nonspher-prop}
Suppose that $\tx{char}(F)\neq 2.$ Let $\alpha$ be as in proposition
\ref{fourier-type}, $\beta\in[0, \frac{\alpha}{4h}),$ and $(V,\tau)$ a non trivial irreducible unitary representation of $K$. There exists $C'>0$, such that the
following holds. Let $C\in\R_+^*$, $(E,\pi)$ any element in
$\mathcal E_{G,C+\beta\ell}$, and $\xi\in E$, $\eta\in V\otimes E^*$ (endowed with the $\ell^2$ norm with respect to some fixed basis of $V$) any
$K$-invariant vectors of norm $1$.  We have for
any $i\geq j\geq 0,$\fixme $$\|c(i,j)\|_V\leq
C'e^{2C-(\frac{\alpha}{2h}-2\beta)i}.$$
\end{prop}

{\bf Proof of theorem \ref{thm-sp4} when $\tx{char}(F)\neq 2$ assuming
proposition
\ref{spher-prop} and \ref{nonspher-prop}:} 
Denote $e_g\in\CC^\E_{C+\beta\ell}(G),\forall g\in G,$ the limit of $\frac{\chi_{E_n}}{vol(E_n)}\in C_c(G)$ for some descending Borel subsets $E_n$ satisfying $\cap_n E_n=\{g\}.$ For any $(\pi,E)\in\E_{G,C+\beta\ell},$ by strong continuity we have $\pi(e_g)\xi=\pi(g)\xi,\forall\xi\in E.$
Let  $P_g=e_Ke_ge_K,$ where $e_K=\int_K e_kdk$ and $dk$ is the Haar measure on
$K$ such that $K$ has volume $1.$ As a consequence of proposition
\ref{spher-prop} we see that the limit $\pr=\lim_{\ell(g)\to\infty}P_g$
exists in $\CC^\E_{C+\beta\ell}(G).$ It is a real and
self-adjoint element because $\bar P_g=P_g,$ and $P_g^*=P_{g^{-1}}.$ Moreover
for any $k\in K$ and $g,g'\in G$ we have $\ell(gkg')\geq\ell(g')-\ell(g^{-1}),$
which gives 
\begin{equation}\label{eqn-spherical}
e_Ke_g\pr=\lim_{\ell(g')\to\infty}e_K\int_K P_{gkg'}dke_K=\pr,
\end{equation}
and therefore $\pr^2=\pr.$

On the other hand, for any non trivial irreducible representation $(V,\tau)$ of
$K$, denote
$e^{V}_K=n\int_K\overline{\text{Tr}(\tau(k))}e_kdk\in\CC_{C+\beta\ell}^\E(G),$
where $n=\dim V.$ For any $(E,\pi)\in\E_{G,C+\beta\ell},$
denote $\pi^*:G\to\mathcal L(E^*)$ the contragredient representation of $\pi,$
i.e. $\pi^*(g)={^t\pi(g^{-1})},$ then 
$\pi^*(e_K^{V})E^*$ is the subspace of vectors in $E^*$ whose $K$-type
is $V.$ For any $\xi\in {\pi^*}(e_K^{V})E^*$ there exist $K$-invariant
vectors $\eta_i\in V^*\otimes E^*$ and vectors $v_i\in V,1\leq i\leq n,$ such
that $\xi=\sum_{i=1}^n\paire{\eta_i}{v_i}.$ By applying proposition
\ref{nonspher-prop} to $V^*$ and $E$ we have $e^V_Ke_ge_K\to 0$ in
$\CC_{C+\beta\ell(g)}^\E$ when $\ell(g)\to\infty,$ and therefore 
\begin{equation}\label{eqn-nonspherical}
e^V_Ke_g\pr=0.
\end{equation}

Note that any vector $z\in E$ satisfying $\pi(e_K^V)z=0$ for any irreducible representation $V$ of $K$
must be the zero vector (since $\pi(f)z=0$ for any class function $f\in C(K)$, i.e. continuous function invariant under the conjugate action of $K$). Now for any $x\in E$ apply this to $z=\pi(e_g\pr-\pr)x,$ and in view of \eqref{eqn-spherical} and \eqref{eqn-nonspherical},
we have
$$\pi(e_g\pr)=\pi(\pr).$$ Therefore $\pi(\pr)E$ is the subspace of
$G$-invariant vectors in $E.$

Finally we complete the proof by taking $\pr_n=P_{D(n,0)}$ and
$t=\frac{\alpha}{2h}-2\beta.$ \qed

Now we turn to the proof of proposition \ref{spher-prop} on spherical matrix
coefficients, which is based on two local estimates on spherical matrix
coefficients corresponding to the move $(0,1)$ and $(1,-1)$ in the Weyl chamber.

\begin{lem}\label{spher-lem-(0,1)}
 Suppose $\tx{char}(F)\neq 2.$ Let $\alpha$ be as in proposition
\ref{fourier-type}. Let $\beta\in[0, \frac{\alpha}{2h})$.
Then there exists $C'>0$, such that for any $C\in\R_+^*$,
any $(E,\pi)\in\mathcal E_{G,C+\beta\ell}$, and  any $K$-invariant vectors 
$\xi\in E$, $\eta\in E^*$ of norm $1$, and any $(i,j)\in\Lambda$ with
$i-j\geq v_0+1,$ we have $$|c(i,j)-c(i,j+1)|\leq
C'e^{2C-(\frac{\alpha}{h}-2\beta) i + \frac{\alpha}{h}j},$$ where $C'$ is a
constant depending on $q,h,v_0,\alpha,\beta.$
\end{lem}

\begin{lem}\label{spher-lem-(1,-1)}
Let $F$ be of any characteristic. Let $\alpha$ be as in proposition
\ref{fourier-type}, and $\beta\in[0, \frac{\alpha}{h})$.
Then there exists $C'>0$, such that for any $C\in\R_+^*$,
any $(E,\pi)\in\mathcal E_{G,C+\beta\ell}$, and  any $K$-invariant vectors 
$\xi\in E$, $\eta\in E^*$ of norm $1$, and $(i,j)\in\Lambda$ with $j\geq
2,$ we have $$|c(i,j)-c(i+1,j-1)|\leq
C'e^{2C+\beta i - (\frac{\alpha}{h}-\beta)j}.$$
\end{lem}

{\bf Proof of proposition \ref{spher-prop} assuming lemma
\ref{spher-lem-(0,1)} and \ref{spher-lem-(1,-1)}:} We adopt the zig-zag argument
from \cite{laf-duke} to $Sp_4.$ We put $$S_\alpha=\{(i,j)\in\Lambda|0\leq
i-2j\leq \alpha\}.$$ First we move any $(i,j)\in\Lambda$ to the strip $S_3.$
Then we show that we can move any $(i,j)\in S_3$ to the line $i=2j$ using the
moves inside $S_4,$ and then we move $(i,j)$ to infinity along this line as
illustrated below.

\vskip 10mm

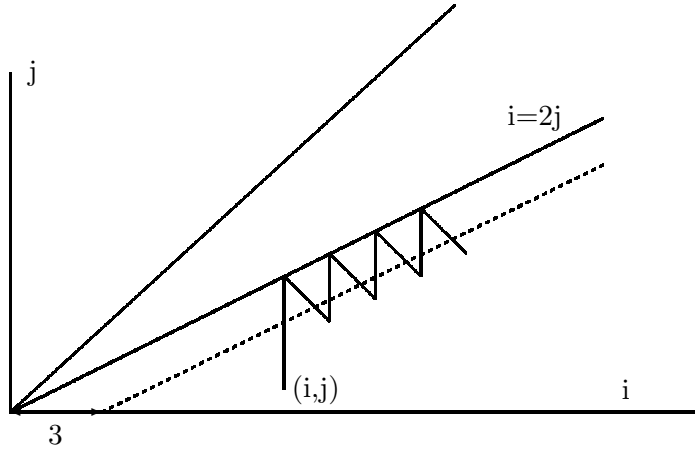
\ifx\JPicScale\undefined\def\JPicScale{0.6}\fi
\unitlength \JPicScale mm
\begin{picture}(155,95)(0,0)
\linethickness{0.3mm}
\put(5,5){\line(1,0){150}}
\linethickness{0.3mm}
\multiput(5,5)(0.13,0.12){750}{\line(1,0){0.13}}
\linethickness{0.3mm}
\put(65,10){\line(0,1){25}}
\linethickness{0.3mm}
\multiput(65,35)(0.12,-0.12){83}{\line(1,0){0.12}}
\linethickness{0.3mm}
\put(75,25){\line(0,1){15}}
\linethickness{0.3mm}
\multiput(75,40)(0.12,-0.12){83}{\line(1,0){0.12}}
\linethickness{0.3mm}
\put(85,30){\line(0,1){15}}
\linethickness{0.3mm}
\multiput(85,45)(0.12,-0.12){83}{\line(1,0){0.12}}
\linethickness{0.3mm}
\put(95,35){\line(0,1){15}}
\linethickness{0.3mm}
\multiput(95,50)(0.12,-0.12){83}{\line(1,0){0.12}}
\put(80,85){\makebox(0,0)[cc]{}}

\put(5,5){\makebox(0,0)[cc]{}}

\linethickness{0.3mm}
\put(5,5){\line(0,1){75}}
\put(10,80){\makebox(0,0)[cc]{j}}

\put(140,10){\makebox(0,0)[cc]{i}}

\put(120,70){\makebox(0,0)[cc]{i=2j}}

\put(72,10){\makebox(0,0)[cc]{(i,j)}}

\linethickness{0.3mm}
\multiput(5,5)(0.24,0.12){542}{\line(1,0){0.24}}
\linethickness{0.3mm}
\multiput(25,5)(1.79,0.89){62}{\multiput(0,0)(0.22,0.11){4}{\line(1,0){0.22}}}
\linethickness{0.3mm}
\put(5,5){\line(1,0){20}}
\put(25,5){\vector(1,0){0.12}}
\put(5,5){\vector(-1,0){0.12}}
\put(15,0){\makebox(0,0)[cc]{3}}

\put(35,-5){\makebox(0,0)[cc]{}}

\end{picture}

\vskip 10mm

Precisely, when $i\geq 2j\geq 0$, we have $(i,\flf{i/2})\in S_2\subset S_3$ and
\begin{align}&|c(i,j)-c(i,\flf{i/2})|\nonumber\\ &\leq C'e^{2C-(\frac{\alpha}{h}-2\beta)i+\frac{\alpha}{h}j} + \ldots
+C'e^{2C-(\frac{\alpha}{h}-2\beta)i+\frac{\alpha}{h}(\flf{i/2}-1)}
\nonumber \\ &\leq
C'e^{2C-(\frac{\alpha}{2h}-2\beta)i}\label{zigzag-1}.\end{align} When
$2j\geq i\geq 0$, we have $\big(i+\clf{\frac{2j-i}{3}},j-\clf{\frac{2j-i}{3}}\big)\in S_3,$ and
\begin{align}&\Big|c(i,j)-c(i+\clf{\frac{2j-i}{3}},j-\clf{\frac{2j-i}{3}})\Big| \nonumber \\&\leq
C'e^{2C-(\frac{\alpha}{h}-\beta)i+\beta
j}+\ldots +C'e^{2C-(\frac{\alpha}{h}-\beta)(i+\clf{\frac{2j-i}{3}}-1)+\beta
(j-\clf{\frac{2j-i}{3}}+1)}\nonumber\\
&\leq C'e^{2C-(\frac{\alpha}{3h}-\beta)(i+j)}\label{zigzag-2}.\end{align}

For any $(i,j)\in S_3,$ if  $i\in2\N+k,k\in\{0,1\}$ then
\begin{align}\big|c(i,j)-c\big(i+k,(i+k)/2\big)\big|\leq
C'e^{2C-(\frac{\alpha}{2h}-2\beta)i}\label{zigzag-3}.\end{align} In fact, by
lemmas \ref{spher-lem-(0,1)} and \ref{spher-lem-(1,-1)}, when $(i,j)\in S_4$ we
have $$\tx{max}\Big(|c(i,j)-c(i,j+1)|,|c(i,j)-c(i+1,j-1)|\Big)\leq
C'e^{2C-(\frac{\alpha}{2h}-2\beta)i}.$$ When $i\in 2\N$ and $(i,j)\in S_3,$ we
get inequality \eqref{zigzag-3} by considering the move $(i,j)\mapsto(i,i/2).$
When $i\in 2\N+1$ and $(i,j)\in S_3,$ there exists $k\in\{0,1\},$ suth that $(i+1,j+k-1)\in
S_4.$ Therefore, we obtain inequality \eqref{zigzag-3} by considering the
following moves inside $S_4:$ $(i,j)\mapsto(i,j+k)\mapsto(i+1,j+k-1)\mapsto(i+1,(i+1)/2).$

Combining inequalities \eqref{zigzag-1}, \eqref{zigzag-2} and \eqref{zigzag-3}
we obtain: when $i\geq 2j\geq 0, $ and $i\in2\N+k,k\in\{0,1\},$
\begin{align}|c(i,j)-c(i+k,(i+k)/2)|\leq
C'e^{2C-(\frac{\alpha}{2h}-2\beta)i}\label{zigzag-4};\end{align} when $2j\geq
i\geq j\geq 0,$ there exists $k\in\{0,1,2\}$ such that
\begin{align}|c(i,j)-c\big(\flf{\frac{2}{3}(i+j)}+k,\frac{1}{2}(\flf{\frac{2}{3}(i+j)}+k)\big)|\leq
C'e^{2C-(\frac{\alpha}{2h}-2\beta)i}\label{zigzag-5}.\end{align}

Finally
for any $j\geq0,$ we have $$|c(2j,j)-c(2j+2,j+1)|\leq
C'e^{2C-(\frac{\alpha}{2h}-2\beta)2j}.$$ Proposition \ref{spher-prop} is then
proved.
\qed

It remains to prove lemmas \ref{spher-lem-(0,1)} and \ref{spher-lem-(1,-1)}. To
prove these two lemmas, we use the following lemma in \cite{laf-jta}
which is a
variant of fast Fourier transform. 
\begin{lem}(lemma 4.4 in \cite{laf-jta})\label{fft}
 Let $\chi:\F\to \C^{*}$ be a non trivial character. Let $h\in \N^{*}, \alpha\in
 \R_{+}^{*}, n\in \N^{*}$. Let $E$ be a Banach space such that
 $\|T_{\resr{h}}\otimes 1_{E}\|\leq e^{-\alpha}$, and let $(\xi_{x,y})_{x,y\in
 \resr{n}}$ be a family of vectors of $E$. Then $$\esp{a,b\in \resr{n}}
 \Big\|\esp{x\in
 \resr{n},\eps \in \F}
 \chi(\eps)\xi_{x,ax+b+\pi^{n-1}\eps}\Big\|^{2} $$ $$\leq q^{2h-2}
 e^{-2(\frac{n}{h}-1)\alpha}\esp{x,y\in \resr{n}} \|\xi_{x,y}\|^{2}.$$
\end{lem}

 {\bf Proof of lemma \ref{spher-lem-(0,1)}:} Denote $m=\flf{\frac{i+j}{2}}$, and
 $n_1=2m-2j-v_0.$ Let $x,y,a,b\in\resr{n_1},$ and let
 $\sigma:\resr{n_1}\to\OO$ be a section. Let $\beta(a,b)^{-1},\alpha(x,y)$ be
 the elements in $G$ defined as follows, $$\beta(a,b)^{-1}= \begin{pmatrix}
 \pi^{m}&&&\\
&\pi^{i-m+j}&& \\
&&\pi^{-i+m-j}& \\
&&&\pi^{-m}
\end{pmatrix}\cdot\begin{pmatrix}
 1&&&\\
0&1&& \\
\sigma(a)&1&1& \\
\sigma(a)^2-2\sigma(b)&\sigma(a)&0&1
\end{pmatrix},$$ $$\alpha(x,y)= \begin{pmatrix} 1&&&\\
0&1&& \\
\sigma(x)&0&1& \\
\sigma(x)^2+2\sigma(y)&\sigma(x)&0&1
\end{pmatrix}\cdot\begin{pmatrix} \pi^{-m+j}&&&\\
&\pi^{-m+j}&& \\
&&\pi^{m-j}& \\
&&&\pi^{m-j}
\end{pmatrix}.$$ Then \begin{gather*}
 \beta(a,b)^{-1}\alpha(x,y)= \begin{pmatrix}
 \pi^{m}&&&\\
&\pi^{i-m+j}&& \\
&&\pi^{-i+m-j}& \\
&&&\pi^{-m}
\end{pmatrix}\times\\\begin{pmatrix}
1&&&\\
0&1&& \\
\sigma(a)+\sigma(x)&1&1& \\
\sigma(a)^2-2\sigma(b)+\sigma(x)^2+2\sigma(y)&\sigma(a)+\sigma(x)&0&1
\end{pmatrix}\cdot\begin{pmatrix} \pi^{-m+j}&&&\\
&\pi^{-m+j}&& \\
&&\pi^{m-j}& \\
&&&\pi^{m-j}
\end{pmatrix}.\end{gather*} 

Recall from the second section that for any
$g\in KD(k,l)K,$ $q^k$ is the biggest norm of all matrix elements in $g,$ and
$q^{k+l}$ is the biggest norm of all $2\times 2$ minors of $g.$
It is easy to see that
$$\|\Lambda^2\big(\beta(a,b)\big)\|=q^{i+j},\|\Lambda^2\big(\alpha(x,y)\big)\|=q^{2m-2j},$$
and $$\|\beta(a,b)^{-1}\alpha(x,y)\|=q^i.$$ On the other hand, we
calculate the minor of rows 3, 4 and columns 1, 2,
\begin{gather*}\det\Big(
\begin{pmatrix}\pi^{-i+m-j}&\\&\pi^{-m}\end{pmatrix}
\begin{pmatrix} \sigma(a)+\sigma(x)&1 \\
\sigma(a)^2-2\sigma(b)+\sigma(x)^2+2\sigma(y)
&\sigma(a)+\sigma(x)
\end{pmatrix}\times\\\begin{pmatrix}\pi^{-m+j}&\\&\pi^{-m+j}\end{pmatrix}\Big)
=-2\pi^{-i-2m+j}\big(\sigma(y)-\sigma(a)\sigma(x)-\sigma(b)\big).
 \end{gather*} Since the norm of the minor of rows 3, 4 and columns 2, 4 is
 $q^{i+j},$ we have $$\|\Lambda^2\big(\beta(a,b)^{-1}\alpha(x,y)\big)\|=\max(q^{i+2m-j-v},q^{i+j}),$$
where $v\in\{0,1,\ldots,2m-2j\}$ is the valuation of
$2(y-ax-b)\in\quotr{2m-2j}.$ Let
$y=ax+b+\pi^{n_1-1}\eps,$ where $\eps\in\F.$ When $\eps=0,$ we have
$v=2m-2j$ and $$\beta(a,b)^{-1}\alpha(x,y)\in KD(i,j)K.$$ When
$\eps\in\F^*$ we have $v=2(m-j)-1$, and then
$$\beta(a,b)^{-1}\alpha(x,y)\in KD(i,j+1)K.$$ 

Let $\chi:\F\to\C^*$ be a non trivial character. By
Cauchy-Schwarz inequality and lemma \ref{fft} we have \begin{align} &|c(i,j)-c(i,j+1)|\nonumber\\
 &=q|\esp{a,b,x\in\resr{n_1},\eps\in\F}
 \chi(\eps)\paire{^t\pi\big(\beta(a,b)\big)\eta}{\pi\big(\alpha(x,ax+b+\pi^{n_1-1}\eps)\big)\xi}|
 \nonumber \\
 &\leq q\sqrt{\esp{a,b\in\resr{n_1}}\|^t\pi\big(\beta(a,b)\big)\eta\|^2}\times
 \nonumber\\
 &\sqrt{\esp{a,b\in\resr{n_1}}\|\esp{x\in\resr{n_1},\eps\in\F}
 \chi(\eps)\pi\big(\alpha(x,ax+b+\pi^{n_1-1}\eps)\big)\xi\|^2} \nonumber \\
 &\leq qe^{C+\beta (i+j)}\cdot q^{h-1}\cdot
 e^{-(\frac{n_1}{h}-1)\alpha}\cdot e^{C+2\beta (m-j)}. \nonumber \\
 &\leq q^{h}\cdot e^{(\frac{v_0+2}{h}+1)\alpha}\cdot
e^{2C-(\frac{\alpha}{h}-2\beta) i + \frac{\alpha}{h}j}, \nonumber
\end{align} and the lemma follows immediately. 
\qed

{\bf Proof of lemma \ref{spher-lem-(1,-1)}:}
Let $x,y,a,b\in\OO/\pi^{j-1} \OO$, and let $\sigma:\resr{j-1}\to\OO$ be a
section.
Define $$\beta(a,b)^{-1}= \begin{pmatrix}
\pi^{i}&&&\\
&1&& \\
&&1& \\
&&&\pi^{-i}
\end{pmatrix}\cdot\begin{pmatrix}
1&&&\\
1+\pi\sigma(a)&1&& \\
0&0&1& \\
-\pi\sigma(b)&0&-1-\pi\sigma(a)&1
\end{pmatrix}\in G,$$
$$\alpha(x,y)= \begin{pmatrix}
1&&&\\
0&1&& \\
\sigma(x)&0&1& \\
\pi\sigma(y)+\sigma(x)&\sigma(x)&0&1
\end{pmatrix}\cdot \begin{pmatrix}
\pi^{-j}&&&\\
&1&& \\
&&1& \\
&&&\pi^{j}
\end{pmatrix}\in G.$$
Then we have 
\begin{gather*}\beta(a,b)^{-1}\alpha(x,y)=\\\begin{pmatrix}
\pi^{i}&&&\\
&1&& \\
&&1& \\
&&&\pi^{-i}
\end{pmatrix}\cdot\begin{pmatrix}
1&&&\\
1+\pi\sigma(a)&1&& \\
\sigma(x)&0&1& \\
\pi(\sigma(y)-\sigma(a)\sigma(x)-\sigma(b)
)&\sigma(x)&-1-\pi\sigma(a)&1
\end{pmatrix}\times\\ \begin{pmatrix}
\pi^{-j}&&&\\
&1&& \\
&&1& \\
&&&\pi^{j}
\end{pmatrix}.\end{gather*}
 
 Firstly, we see that
 $$\|\Lambda^2\big(\beta(a,b)\big)\|=q^{i},\|\Lambda^2\big(\alpha(x,y)\big)\|=q^{j},$$
 and $$\|\Lambda^2\big(\beta(a,b)^{-1}\alpha(x,y)\big)\|=q^{i+j},$$ which is the
 norm of the determinant of the submatrix of rows 2, 4 and columns 1, 3.
 Denote the valuation of $y-ax-b\in\quotr{j-1}$ by $v\in\{0,1,\ldots,j-1\},$
 and we have
 $$\|\beta(a,b)^{-1}\alpha(x,y)\|=\max(q^i,q^{i+j-v-1}).$$ Let
 $$y=ax+b+\pi^{j-2}\eps,\eps\in \F^.$$ When $\eps=0,$ we see that $v=j-1$ and
 $$\beta(a,b)^{-1}\alpha(x,ax+b)\in KD(i,j)K.$$ When $\eps\in\F^*$ we
 have $v=j-2$, and therefore $$\beta(a,b)^{-1}\alpha(x,ax+b+\pi^{j-2}\eps)\in
 KD(i+1,j-1)K.$$
 
 Let $\chi:\F\to \C^*$ be a non trivial character. By the same estimates as in
 the end of the proof of lemma \ref{spher-lem-(0,1)} ($n_1$ replaced by $j-1$),
 we have
 \begin{align}
 &|c(i,j)-c(i+1,j-1)| \nonumber\\ \leq& qe^{C+\beta j}\cdot q^{h-1}\cdot
 e^{-(\frac{j-1}{h}-1)\alpha}\cdot
  e^{C+\beta i} \nonumber \\ =&q^{h}\cdot e^{(\frac{1}{h}+1)\alpha}\cdot
  e^{2C+\beta i - (\frac{\alpha}{h}-\beta)j}. \nonumber \end{align} 
\qed

As for proposition \ref{nonspher-prop}, we need two similar lemmas as follows
for its proof. 

\begin{lem}\label{nonspher-lem-(0,1)}
Suppose $\tx{char}(F)\neq 2.$  Let $\alpha$ be as in proposition
\ref{fourier-type}, $\beta\in[0, \frac{\alpha}{2h}),$ and $(V,\tau)$ a non
trivial irreducible unitary representation of $K$ which factorizes through
$Sp_4(\resr{k})$ for $k\geq 1$. There exists $C'>0$, such that the following
holds. Let $C\in\R_+^*$, $(E,\pi)$ any element in $\mathcal E_{G,C+\beta\ell}$, and $\xi\in E$, $\eta\in V\otimes E^*$ any $K$-invariant vectors of norm $1$.  Then
for any $(i,j)\in\Lambda$ with $i- j \geq 2k+v_0,$ we have
$$\|c(i,j)-c(i,j+1)\|_V\leq
C'e^{2C-(\frac{\alpha}{h}-2\beta) i + \frac{\alpha}{h}j}.$$
\end{lem}

\begin{lem}\label{nonspher-lem-(1,-1)}
Let $F$ be of any characteristic. Let $\alpha$ be as in proposition
\ref{fourier-type}, $\beta\in[0, \frac{\alpha}{h}),$ and $(V,\tau)$ a non
trivial irreducible unitary representation of $K$ which factorizes through
$Sp_4(\resr{k})$ for $k\geq 1$. There exists $C'>0$, such that the following
holds. Let $C\in\R_+^*$, $(E,\pi)$ any element in $\mathcal E_{G,C+\beta\ell}$, and $\xi\in E$, $\eta\in V\otimes E^*$ any
$K$-invariant vectors of norm $1$.  Then for any $(i,j)\in\Z^2$
with $i+1\geq j\geq 2k+2,$ we have
$$\|c(i,j)-c(i+1,j-1)\|_V\leq
C'e^{2C+\beta i - (\frac{\alpha}{h}-\beta)j}.$$ In particular,
$$\|c(j-1,j)-c(j,j-1)\|_V\leq C'e^{2C- (\frac{\alpha}{h}-2\beta)j}.$$ 
\end{lem}


\begin{lem}\label{fft-with-k}
Let $h,\alpha,n,E$ as in lemma \ref{fft}. Let $k\in \{0,...,\flf{n/2}\}$,
$\eps_0\in\F^*$, and let $(\xi_{x,y})_{x\in\pi^k\resr{n},y\in\pi^{2k}
 \resr{n}}$ be a family of vectors of $E$. Then there exists a constant $C_2$
 depending only on $q,$ such that $$\esp{a\in\pi^k\resr{n},b\in
 \pi^{2k}\resr{n}} \Big\|\esp{x\in \pi^{k}\resr{n}}\xi_{x,ax+b+\pi^{n-1}\eps_0}-\esp{x\in
 \pi^{k}\resr{n}}\xi_{x,ax+b}  \Big\|^{2}$$ $$\leq
 C_2q^{2h-2}e^{-2(\frac{n-2k}{h}-1)\alpha}\esp{x\in
 \pi^{k}\resr{n},y\in\pi^{2k}\resr{n}}\|\xi_{x,y}\|^{2}.$$   
\end{lem}

{\bf Proof:} When $k=0$, let $f$ be the function on $\F$ defined by
$f(\eps_0)=q,$ $f(0)=-q,$ and zero elsewhere. The left hand side of the inequality is equal to
$$\esp{a,b\in\quotr{n}}\|\esp{x\in\quotr{n},\eps\in\F}f(\eps)\xi_{x,ax+b+\pi^{n-1}\eps}\|^2.$$
Write $f=\sum_{\chi\in\widehat{\F},\chi\neq 1}f_\chi\chi $ with
$f_\chi\in\C,$ then by the triangular inequality and and lemma \ref{fft}, the left hand side is equal to
$$\esp{a,b\in\quotr{n}}\|\sum_{\chi\in\widehat{\F},\chi\neq
1}f_\chi\esp{x\in\quotr{n},\eps\in\F}\chi(\eps)\xi_{x,ax+b+\pi^{n-1}\eps}\|^2$$
$$\leq
C_2\max_{\chi\in\widehat{\F},\chi\neq
1}\esp{a,b\in\quotr{n}}\|\esp{x\in\quotr{n},\eps\in\F}\chi(\eps)\xi_{a,ax+b+\pi^{n-1}\eps}\|^2$$
$$\leq C_2
q^{2h-2}e^{-2(\frac{n}{h}-1)\alpha}\esp{a,b\in\quotr{n}}\|\xi_{x,y}\|^2,$$ where
$C_2=\big(\sum_{\chi\in\widehat{\F},\chi\neq 1}|f_\chi|\big)^2.$ 

In general, let $s:\resr{n-2k}\to\resr{n-k}$ be a section, and for any
$x_1,y_1\in\quotr{n-2k}$ let $$\xi'_{x_1,y_1}=\esp{z\in\pi^{n-2k}\quotr{n-k}}\xi_{\pi^{k}(s(x_1)+z),\pi^{2k}y_1}.$$ For any $a,x\in\pi^k\quotr{n},$ the product $ax\in\quotr{n}$ only depends on
 the images of $a,x$ in $\pi^k\quotr{n-k}$. So the left hand side of the
 inequality is equal to
 $$\esp{a_1,b_1\in\quotr{n-2k}}\|\esp{x_1\in\quotr{n-2k}}\big(\xi'_{x_1,a_1x_1+b_1+\pi^{n-2k-1}\eps_0}-\xi'_{x_1,a_1x_1+b_1}\big)\|^2.$$
 By applying the lemma when $k=0$ to $(\xi'_{x_1,y_1})_{x_1,y_1\in\quotr{n-2k}}$
 we get the inequality in the lemma with the same $C_2.$
  \qed  

{\bf Proof of lemma \ref{nonspher-lem-(0,1)}:}

Let $m,x,y,a,b,\eps,\sigma,\alpha(x,y),\beta(a,b)$ be as in the proof of
lemma \ref{spher-lem-(0,1)}. We recall also from the proof that
$$\|\Lambda^2\big(\beta(a,b)\big)\|=q^{i+j},\|\Lambda^2\big(\alpha(x,y)\big)\|=q^{2m-2j}.$$

Let $\eps_0$ be image of $\pi^{v_0}/2$ in $\F^*$, and let
$$\varepsilon_1=2\pi^{-2m+2j+1}(\sigma(y)-\sigma(a)\sigma(x)-\sigma(b))\in\OO.$$
Recall that $y=ax+b+\pi^{2m-2j-v_0-1}\eps$, we have
$\varepsilon_1\tx{ mod}\pi\OO=\eps_0^{-1}\varepsilon .$ Let $k_1$ be the element
in $K$ defined by $$k_1=\begin{pmatrix} 0&0&1&0\\
0&0&-\pi^{i-2m+j}(\sigma(a)+\sigma(x))&1
\\
-1&0&-\pi^{i-2m+3j+1}(\sigma(a)+\sigma(x))&\pi^{2j+1} \\
-\pi^{i-2m+j}(\sigma(a)+\sigma(x))&-1&\pi^{2i-2m+2j}&0 \end{pmatrix},$$ and let
$$g_1=k_1\beta(a,b)^{-1}\alpha(x,y)$$ $$=\begin{pmatrix} 0&0&1&0 \\
0&0&-\pi^{i-2m+j}(\sigma(a)+\sigma(x))&1 \\
-1&0&-\pi^{i-2m+3j+1}(\sigma(a)+\sigma(x))&\pi^{2j+1} \\
-\pi^{i-2m+j}(\sigma(a)+\sigma(x))&-1&\pi^{2i-2m+2j}&0 \end{pmatrix} \times$$ $$ \begin{pmatrix} \pi^{j}&&&\\
0&\pi^{i-2m+2j}&& \\
\pi^{-i}(\sigma(a)+\sigma(x))&\pi^{-i}&\pi^{-i+2m-2j}& \\
\pi^{-2m+j}\big(\sigma(a)+\sigma(x)\big)^2+\pi^{-j-1}\eps_1&\pi^{-2m+j}(\sigma(a)+\sigma(x))&0&\pi^{-j}
\end{pmatrix}$$ $$=\begin{pmatrix}
\pi^{-i}(\sigma(a)+\sigma(x))&\pi^{-i}&\pi^{-i+2m-2j}&0\\
\pi^{-j-1}\eps_1&0&-\pi^{-j}(\sigma(a)+\sigma(x))&\pi^{-j} \\
\pi^{j}(\eps_1-1)&0&-\pi^{j+1}(\sigma(a)+\sigma(x))&\pi^{j+1}\\
0&0&\pi^{i}&0 \end{pmatrix}.$$ When $\eps=0,$ we have
$$\begin{pmatrix} \pi^i&&&\\ &\pi^{j}&& \\ &&\pi^{-j}& \\
&&&\pi^{-i}\end{pmatrix} g_1=\begin{pmatrix}
\sigma(a)+\sigma(x)&1&\pi^{2m-2j}&0\\
\pi^{-1}\eps_1&0&-(\sigma(a)+\sigma(x))&1 \\
\eps_1-1&0&-\pi(\sigma(a)+\sigma(x))&\pi\\
0&0&1&0 \end{pmatrix}$$ which is an element in $K.$ When $\eps=\eps_0,$ we have
$$\begin{pmatrix}\pi^i&&&\\ &\pi^{j+1}&& \\ &&\pi^{-j-1}& \\
&&&\pi^{-i}\end{pmatrix} g_1=\begin{pmatrix}
\sigma(a)+\sigma(x)&1&\pi^{2m-2j}&0\\
\eps_1&0&-\pi(\sigma(a)+\sigma(x))&\pi \\
\pi^{-1}(\eps_1-1)&0&-(\sigma(a)+\sigma(x))&1\\
0&0&1&0 \end{pmatrix},$$ which is also in $K.$ Denote
$\xi_{x,y}=\pi(\alpha(x,y))\xi,$
$\eta_{a,b}=\big(\text{Id}_V\otimes^t\pi(\beta(a,b))\big)\eta$ and
$n_1=2(m-j)-v_0.$ Note that $c(k'gk'')=\tau(k')c(g)$ for any $k',k''\in K,g\in
G$, we then have \begin{gather} \|c(i,j)-c(i,j+1)\|_V\nonumber\\
 =q\Big\|\esp{a,x\in\pi^k\resr{n_1},b\in\pi^{2k}\resr{n_1}}
 \tau(k_1)\Big(\paire{\eta_{a,b}}{\xi_{x,ax+b+\pi^{n_1-1}\eps_0}}-\paire{\eta_{a,b}}{\xi_{x,ax+b}}\Big)\Big\|_V
 .\label{norm}\end{gather}
 
 When $i-j\geq k+v_0$ and $a,x\in\pi^k\resr{n_1},$
 we have
 $$k_1=\begin{pmatrix}0&0&1&0\\0&0&0&1\\-1&0&0&\pi^{2j+1}\\0&-1&0&0\end{pmatrix}\tx{
 mod}\pi^k\OO,$$ so (\ref{norm}) becomes $$ q
 \Big\|\esp{a,x\in\pi^k\resr{n_1},b\in\pi^{2k}\resr{n_1}} \Big(\paire{\eta_{a,b}}{\xi_{x,ax+b+\pi^{n_1-1}\eps_0}}-\paire{\eta_{a,b}}{\xi_{x,ax+b}}\Big)\Big\|_V.$$
 By Cauchy-Schwarz inequality and lemma \ref{fft-with-k} (when $i-j\geq
 2k+v_0$), it is less than 
\begin{align}
 &\leq
 q\sqrt{\esp{a\in\pi^k\resr{n_1},b\in\pi^{2k}\resr{n_1}}\|\eta_{a,b}\|^2}\times \nonumber\\
 &\sqrt{\esp{a\in\pi^k\resr{n_1},b\in\pi^{2k}\resr{n_1}}\|\esp{x\in\pi^k\resr{n_1}}
 \xi_{x,ax+b+\pi^{n_1-1}\eps_0-}\esp{x\in\pi^k\resr{n_1}}
 \xi_{x,ax+b}\|^2} \nonumber \\
 &\leq qe^{C+\beta (i+j)}\cdot \sqrt{C_2}q^{h-1}\cdot
 e^{-(\frac{n_1-2k}{h}-1)\alpha}\cdot e^{C+2\beta (m-j)} \nonumber \\
&\leq \sqrt{C_2}q^{h}\cdot e^{(\frac{v_0+2+2k}{h}+1)\alpha}\cdot
e^{2C-(\frac{\alpha}{h}-2\beta) i + \frac{\alpha}{h}j}. \nonumber
\end{align}
\qed 

{\bf Proof of lemma \ref{nonspher-lem-(1,-1)}:} 

Let
$x,y,a,b,\eps,\sigma,\alpha(x,y),\beta(a,b)$ be as in the proof of lemma
\ref{spher-lem-(1,-1)}. From the proof we have
$$\|\Lambda^2\big(\beta(a,b)\big)\|=q^{i},\|\Lambda^2\big(\alpha(x,y)\big)\|=q^{j}.$$

Denote
$\varepsilon_1=\pi^{-j+2}(\sigma(y)-\sigma(a)\sigma(x)-\sigma(b))\in\OO,$
and we have $\eps_1\tx{ mod }\pi\OO=\eps.$ Denote $a_1=1+\pi\sigma(a)\in\OO.$
For any $i+1\geq j\geq 1,$ let $k_1$ be the element in $K$ defined by
$$k_1=\begin{pmatrix}
0&0&0&1
\\
0&1&0&-\pi^{i-j+1}a_1 \\
0&-\pi^{j-1}a_1^{-2}\eps_1-a_1^{-1}\sigma(x)&1&\pi^ia_1^{-1} \\
-1&\pi^ia_1^{-1}(1-\eps_1)-\pi^{i-j+1}\sigma(x)&\pi^{i-j+1}a_1&\pi^{2i-j+1}
\end{pmatrix}.$$ Denote $$g_1=k_1\beta(a,b)^{-1}\alpha(x,y)$$ $$=\begin{pmatrix}
0&0&0&1
\\
0&1&0&-\pi^{i-j+1}a_1 \\
0&-\pi^{j-1}a_1^{-2}\eps_1-a_1^{-1}\sigma(x)&1&\pi^ia_1^{-1} \\
-1&\pi^ia_1^{-1}(1-\eps_1)-\pi^{i-j+1}\sigma(x)&\pi^{i-j+1}a_1&\pi^{2i-j+1}
\end{pmatrix} \times$$ $$\begin{pmatrix} \pi^{i-j}&&&\\
\pi^{-j}a_1&1&& \\
\pi^{-j}\sigma(x)&0&1& \\
\pi^{-i-1}\eps_1&\pi^{-i}\sigma(x)&-\pi^{-i}a_1&\pi^{-i+j}
\end{pmatrix}$$ $$=\begin{pmatrix}
\pi^{-i-1}\eps_1&\pi^{-i}\sigma(x)&-\pi^{-i}a_1&\pi^{-i+j} \\
\pi^{-j}a_1(1-\eps_1)&1-\pi^{-j+1}a_1\sigma(x)
&\pi^{-j+1}a_1^2&-\pi a_1 \\ 
0&-\pi^{j-1}a_1^{-2}\eps_1&0&\pi^ja_1^{-1} \\
0&\pi^{i}a_1^{-1}(1-\eps_1)&0&\pi^{i+1} \end{pmatrix}$$ When
$\eps=0,$ we have $$\begin{pmatrix}\pi^i&&&\\ &\pi^{j}&& \\ &&\pi^{-j}& \\
&&&\pi^{-i}\end{pmatrix} g_1$$ $$= \begin{pmatrix}
\pi^{-1}\eps_1&\sigma(x)&-a_1&\pi^{j} \\
a_1(1-\eps_1)&\pi^j-\pi a_1\sigma(x)
&\pi a_1^2&-\pi^{j+1}a_1 \\ 
0&-\pi^{-1}a_1^{-2}\eps_1&0&a_1^{-1} \\
0&a_1^{-1}(1-\eps_1)&0&\pi \end{pmatrix}\in K.$$ When
$\eps=1,$ we have $$\begin{pmatrix}\pi^{i+1}&&&\\ &\pi^{j-1}&& \\
&&\pi^{-j+1}& \\
&&&\pi^{-i-1}\end{pmatrix} g_1$$ $$= \begin{pmatrix}
\eps_1&\pi\sigma(x)&-\pi a_1&\pi^{j+1} \\
\pi^{-1}a_1(1-\eps_1)&\pi^{j-1}-a_1\sigma(x)
&a_1^2&-\pi^{j}a_1 \\ 
0&-a_1^{-2}\eps_1&0&\pi a_1^{-1} \\
0&\pi^{-1}a_1^{-1}(1-\eps_1)&0&1 \end{pmatrix}\in K$$ When $j\geq 2k+2$ and
$a,x\in\pi^k\resr{j-1},$ we have
$$k_1=\begin{pmatrix}0&0&0&1\\0&1&0&-\pi^{i-j+1}\\0&0&1&0\\-1&0&\pi^{i-j+1}&0\end{pmatrix}\tx{
mod }\pi^k\OO.$$ By the same estimates (with $n_1$ replaced by $j-1$ and
$\eps_0$ by $1\in\F^*$) as in the end of the proof of lemma
\ref{nonspher-lem-(0,1)} we have $$\|c(i,j)-c(i+1,j-1)\|_V\leq qe^{C+\beta j}\cdot \sqrt{C_2}q^{h-1}\cdot
e^{-(\frac{j-1-2k}{h}-1)\alpha}\cdot e^{C+\beta i}$$ $$= \sqrt{C_2}q^{h}\cdot
e^{(\frac{1+2k}{h}+1)\alpha}\cdot e^{2C+\beta i-(\frac{\alpha}{h}-\beta)j}.$$
\qed 

Let $K_1$ be the subgroup of $K$ consisting of elements of the form
$$\begin{pmatrix}*&*&& \\ *&*&& \\&&*&*\\&&*&*\end{pmatrix}$$, and $K_2$
consisting of elements of the form
$$\begin{pmatrix}1&&&\\&*&*&\\&*&*&\\&&&1\end{pmatrix},$$ i.e.
$$K_1=\{\begin{pmatrix}A&\\&
Q^tA^{-1}Q\end{pmatrix}|A\in GL_2(\OO)\},$$ where
$Q=\begin{pmatrix}&1\\1&\end{pmatrix},$ and
$$K_2=\{\begin{pmatrix}1&&\\&B&\\&&1\end{pmatrix}|B\in SL_2(\OO)\}.$$ 

\begin{lem}\label{generated-sg}
Let $F$ be of any characteristic. Then $K=(K_1K_2)^{30}.$
\end{lem}

{\bf Proof:} Denote $B$ the lower triangular matrices in $K,$ and $W$ the Weyl
group associated to $G=Sp_4(F).$ Denote
$$w_{21}=\begin{pmatrix}0&1&&\\1&0&&\\&&0&1\\&&1&0\end{pmatrix},
w_{32}=\begin{pmatrix}1&&&\\&0&1&\\&-1&0&\\&&&1\end{pmatrix}.$$  The dihedral
group $W$ (of order $8$) is generated by $w_{21}$ and $w_{32},$ which are
reflections w.r.t. the axes $x=y$ and $x=0,$ respectively. Since $w_{21}\in K_1$
and $w_{32}\in K_2$ we obtain $W\subset (K_1K_2)^{4}.$ 

Denote for any $a\in \OO,$
\begin{gather*}\mu_{21}(a)=\begin{pmatrix}1&0&&\\a&1&&\\&&1&0\\&&-a&1\end{pmatrix}
,\mu_{32}(a)=\begin{pmatrix}1&&&\\&1&0&\\&a&1&\\&&&1\end{pmatrix},\\
\mu_{31}(a)=\begin{pmatrix}1&&&\\0&1&&\\a&0&1&\\0&a&0&1\end{pmatrix},
\mu_{41}(a)=\begin{pmatrix}1&&&\\0&1&&\\0&0&1&\\a&0&0&1\end{pmatrix}.\end{gather*}
By calculations we have
$$\mu_{41}(a)=w_{21}\mu_{32}(a)w_{21}\in(K_1K_2)^3$$ and
$$\mu_{31}(a)=\mu_{21}(-a)\mu_{32}(1)\mu_{21}(a)\mu_{32}(-1)\mu_{41}(-a^2)\in(K_1K_2)^7.$$ Any element in $B$ has the form
$$\begin{pmatrix}1&&&\\a&1&&\\c&b&1&\\d&c-ab&-a&1\end{pmatrix}\cdot
\begin{pmatrix}e&&&\\&f&&\\&&f^{-1}&\\&&&e^{-1}\end{pmatrix}$$ where
$a,b,c,d\in\OO$ and $e,f\in\OO^\times,$ which is equal to
$$\mu_{21}(a)\mu_{32}(b)\mu_{31}(c)\mu_{41}(ac+d)\cdot
\begin{pmatrix}e&&&\\&f&&\\&&f^{-1}&\\&&&e^{-1}\end{pmatrix}.$$ So we have
$B\subset (K_1K_2)^{13}.$

By the Bruhat decomposition, we have $K=BWB=(K_1K_2)^{30}.$\qed

\begin{lem}\label{invariant}
Let $K$ be any compact group, $\{K_i\}_{1\leq i\leq
n}$ a family of subgroups such that $K=(K_1K_2\ldots K_n)^N$ for some
$N\in\N^*.$ Then for any finite dimensional unitary representation $(V,\tau)$ of
$K$ without invariant vector, and any $x\in V$, and $y_i\in V$ invariant by
$K_{i}$ for each $1\leq i\leq n,$ we have $$\|x\|_{V}\leq 2nN\max_{1\leq i\leq
n}\{\|x-y_i\|_{V}\}.$$  
\end{lem}
{\bf Proof:} Since $\int_K\|\tau(k)x-x\|_{V}^2dk=2\|x\|_{V}^2\geq\|x\|_{V}^2$ we
see that there exists a $k\in K$ such that $\|\tau(k)x-x\|_{V}\geq\|x\|_{V}.$
Suppose that $k=(k_{11}\ldots k_{n1})\ldots(k_{1N}\ldots k_{nN})$ with $k_{ij}\in K_i(1\leq i\leq
n,1\leq j\leq N).$ We then have $$\|x\|_{V}\leq\|\tau(k)x-x\|_{V}\leq
\sum_{1\leq i\leq n,1\leq j\leq N}\|\tau(k_{ij})x-x\|_{V}$$$$\leq 2\sum_{1\leq
i\leq n,1\leq j\leq N}\|y_i-x\|_{V}\leq 2nN\max_{1\leq i\leq
n}\{\|x-y_i\|_{V}\}$$ 
\qed

{\bf Proof of proposition \ref{nonspher-prop}:} 
By lemmas \ref{nonspher-lem-(0,1)} and \ref{nonspher-lem-(1,-1)}, we obtain two
similar inequalities as \eqref{zigzag-4} and \eqref{zigzag-5} in the proof of
proposition \ref{spher-prop} (using the same argument): when $i\geq 2j\geq 0, $ and $i\in2\N+k,k\in\{0,1\},$
$$\|c(i,j)-c(i+k,(i+k)/2)\|_V\leq
C'e^{2C-(\frac{\alpha}{2h}-2\beta)i};$$ when $2j\geq
i\geq j\geq 0,$ there exists $k\in\{0,1,2\}$ such that
$$\|c(i,j)-c\big(\flf{\frac{2}{3}(i+j)}+k,\frac{1}{2}(\flf{\frac{2}{3}(i+j)}+k)\big)\|_V\leq
C'e^{2C-(\frac{\alpha}{2h}-2\beta)i}.$$ So it remains to prove
$$\|c(2j,j)\|_V\leq C'e^{2C-\frac{\alpha}{h}-2\beta)2j}.$$
First we see that
 $$\tx{max}\Big(\|c(2j,j)-c(2j,0)\|_V,\|c(2j,j)-c(\flf{3j/2},\flf{3j/2})\|_V\Big)$$$$\leq
 C'e^{2C-(\frac{\alpha}{2h}-2\beta)2j}. $$ Moreover,
 $c(k'gk'')=\tau(k')c(g),\forall k',k''\in K,g\in G,$ and it follows that 
 $c(\flf{3j/2},\flf{3j/2})$ is invariant by $K_1,$ and that $c(2j,0)$ invariant
 by $K_2.$ By applying lemma \ref{invariant} to $K=(K_1K_2)^{30}$, we complete
 the proof of the proposition.\qed

\section{Proof of theorem \ref{thm-sp4} when $\tx{char}(F)=2$}
In this section we prove theorem \ref{thm-sp4} when $\tx{char}(F)=2.$ The proof
for $\tx{char}(F)=2$ is technically more difficult because it is only possible
to prove a local estimate for the move $(0,2)$, and
therefore we have two limits in the spherical propositions (proposition
\ref{spher-prop-char2}).

Throughout this section we assume $F$ is of characteristic $2.$

\begin{lem}\label{ch2-lem-(0,2)} Let $\alpha>0$ as in proposition
\ref{fourier-type}, $\beta\in[0, \frac{\alpha}{4h}).$ Let $(V,\tau)$ be an irreducible unitary
representation of $K$ which factorizes through $Sp_4(\resr{k})$ for $k\geq 0$.
There exists $C'>0$, such that the following holds.
Let $C\in\R_+^*$, $(E,\pi)$ any element in $\mathcal E_{G,C+\beta\ell}$, and
$\xi\in E$, $\eta\in V\otimes E^*$ any $K$-invariant vectors of norm $1$. Then
for any $(i,j)\in\Lambda$ with $i- j \geq 4k+2,$ we have
$$\|c(i,j)-c(i,j+2)\|_V\leq
C'e^{2C-(\frac{\alpha}{2h}-2\beta) i + \frac{\alpha}{2h}j}.$$ In particular when
$(V,\tau)$ is the trivial representation of $K$ (and $V=\C$), we have
$$|c(i,j)-c(i,j+2)|\leq C'e^{2C-(\frac{\alpha}{2h}-2\beta) i +
\frac{\alpha}{2h}j},$$ for any $(i,j)\in\Lambda$ with $i-j\geq 2.$ 
\end{lem}

{\bf Proof:} Since $\text{char}(F)=2,$ we have $-1=1$ in $F.$ Let
$m=\flf{\frac{i+j}{2}}$, $x,y,a,b\in\resr{m-j-1}$ satisfying $y+ax+b\in\pi^{m-j-2}\resr{m-j-1}$, and put $\eps\in\F$ with $y=ax+b+\pi^{m-j-2}\eps$ as usual. Let 
$\sigma:\resr{m-j-1}\to\OO$ be a
section. Let $$\beta(a,b)^{-1}= \begin{pmatrix} \pi^{m}&&&\\
0&\pi^{i-m+j}&& \\
\pi^{-i+m-j+1}\sigma(b)&\pi^{-i+m-j}(1+\pi \sigma(a))^2&\pi^{-i+m-j}& \\
0&\pi^{-m+1}\sigma(b)&0&\pi^{-m} \end{pmatrix},$$ $$\alpha(x,y)=
\begin{pmatrix} \pi^{-m+j}&&&\\
0&\pi^{-m+j}&& \\
\pi^{-m+j}(\sigma(x)+\pi \sigma(y))&0&\pi^{m-j}& \\
\pi^{-m+j}\sigma(x)^2&\pi^{-m+j}(\sigma(x)+\pi \sigma(y))&0&\pi^{m-j}
\end{pmatrix}.$$ Then $$\beta(a,b)^{-1}\alpha(x,y)=$$ $$ \begin{pmatrix}
\pi^{j}&&&\\
0&\pi^{i-2m+2j}&& \\
\pi^{-i}(\pi\sigma(b)+\sigma(x)+\pi \sigma(y))&\pi^{-i}(1+\pi
\sigma(a))^2&\pi^{-i+2m-2j}& \\
\pi^{-2m+j}\sigma(x)^2&\pi^{-2m+j}(\pi\sigma(b)+\sigma(x)+\pi
\sigma(y))&0&\pi^{-j} \end{pmatrix}.$$ 

 We see that
 $$\|\Lambda^2\big(\beta(a,b)\big)\|=q^{i+j},\|\Lambda^2\big(\alpha(x,y)\big)\|=q^{2m-2j}.$$
 Denote $$a_1=(\pi\sigma(b)+\sigma(x)+\pi\sigma(y))(1+\pi\sigma(a))^{-2},$$
 and $$\varepsilon_1
 =\pi^{-m+j+2}\big(\sigma(y)+\sigma(a)\sigma(x)+\sigma(b)\big)\in\OO.$$ Let
 $k_1$ be the element in $K$ defined by $$k_1=\begin{pmatrix} 0&0&1&0\\
0&0&\pi^{i-2m+j}a_1&1 \\
1&0&\pi^{i-2m+3j+2}a_1&\pi^{2j+2} \\
\pi^{i-2m+j}a_1&1&\pi^{2i-2m+2j}(1+\pi\sigma(a))^{-2}&0 \end{pmatrix},$$ and let
$$g_1=k_1\beta(a,b)^{-1}\alpha(x,y)$$ $$=\begin{pmatrix} 0&0&1&0 \\
0&0&\pi^{i-2m+j}a_1&1 \\
1&0&\pi^{i-2m+3j+2}a_1&\pi^{2j+2} \\
\pi^{i-2m+j}a_1&1&\pi^{2i-2m+2j}(1+\pi\sigma(a))^{-2}&0 \end{pmatrix} \times$$
$$ \begin{pmatrix} \pi^{j}&&&\\
0&\pi^{i-2m+2j}&& \\
\pi^{-i}(\pi\sigma(b)+\sigma(x)+\pi \sigma(y))&\pi^{-i}(1+\pi
\sigma(a))^2&\pi^{-i+2m-2j}& \\
\pi^{-2m+j}\sigma(x)^2&\pi^{-2m+j}(\pi\sigma(b)+\sigma(x)+\pi \sigma(y))&0&\pi^{-j}
 \end{pmatrix}$$ $$=\begin{pmatrix}
\pi^{-i}a_1(1+\pi\sigma(a))^{2}&\pi^{-i}(1+\pi\sigma(a))^{2}&\pi^{-i+2m-2j}&0\\
\pi^{-j-2}\eps_1^2(1+\pi\sigma(a))^{-2}&0&\pi^{-j}a_1&\pi^{-j} \\
\pi^{j}\eps_1^2(1+\pi\sigma(a))^{-2}+\pi^j&0&\pi^{j+2}a_1&\pi^{j+2}\\
0&0&\pi^{i}(1+\pi\sigma(a))^{-2}&0 \end{pmatrix}.$$ When $\eps=0,$ we have
$|\eps_1^2|\leq q^{-2}$ and $$\begin{pmatrix} \pi^i&&&\\ &\pi^{j}&& \\
&&\pi^{-j}&
\\
&&&\pi^{-i}\end{pmatrix} g_1$$$$=\begin{pmatrix}
a_1(1+\pi\sigma(a))^{2}&(1+\pi\sigma(a))^{2}&\pi^{2m-2j}&0\\
\pi^{-2}\eps_1^2(1+\pi\sigma(a))^{-2}&0&a_1&1 \\
\eps_1^2(1+\pi\sigma(a))^{-2}+1&0&\pi^{2}a_1&\pi^{2}\\
0&0&(1+\pi\sigma(a))^{-2}&0 \end{pmatrix}\in K.$$ When $\eps=1,$ we have
$$|\eps_1^2(1+\pi\sigma(a))^{-2}+1|
=\big|\big(\eps_1(1+\pi\sigma(a))^{-1}+1\big)^2\big|\leq q^{-2},$$ and then
$$\begin{pmatrix}\pi^i&&&\\ &\pi^{j+2}&& \\ &&\pi^{-j-2}& \\  
&&&\pi^{-i}\end{pmatrix} g_1$$ $$=\begin{pmatrix}
a_1(1+\pi\sigma(a))^{2}&(1+\pi\sigma(a))^{2}&\pi^{2m-2j}&0\\
\eps_1^2(1+\pi\sigma(a))^{-2}&0&\pi^{2}a_1&\pi^{2} \\
\pi^{-2}\big(\eps_1^2(1+\pi\sigma(a))^{-2}+1\big)&0&a_1&1\\
0&0&(1+\pi\sigma(a))^{-2}&0 \end{pmatrix}\in K.$$ When $i-j\geq 4k+2$
, $a,x\in\pi^k\resr{m-j-1},$ and $b,y\in\pi^{2k}\resr{m-j-1}$, we have
$$k_1=\begin{pmatrix} 0&0&1&0\\
0&0&0&1 \\
1&0&0&\pi^{2j+2} \\
0&1&0&0 \end{pmatrix} \tx{ mod }\pi^{k}\OO.$$ By replacing $n_1$ by
$m-j-1$ at the end of the proof of lemma \ref{nonspher-lem-(0,1)} we get
\begin{gather*}\|c(i,j)-c(i,j+2)\|_V\\ \leq qe^{C+\beta(i+j)}\cdot
\sqrt{C_2}q^{h-1}\cdot e^{-(\frac{m-j-1-2k}{h}-1)\alpha}\cdot e^{C+2\beta(m-j)}\\
\leq C'e^{2C-(\frac{\alpha}{2h}-2\beta)i+\frac{\alpha}{2h}j}.\end{gather*}
\qed

\begin{prop}\label{spher-prop-char2}
Let $\alpha>0$, $\beta\in[0, \frac{\alpha}{8h})$. There exists $C'>0$, such
that the following holds. Let $C\in\R_+^*$, $(E,\pi)$ any
element in $\mathcal E_{G,C+\beta\ell}$, and $\xi\in E$, $\eta\in E^*$ any
$K$-invariant vectors of norm $1$. There exist $c_0,c_1\in\C$, such that $$|c(i,j)-c_l|\leq
C'e^{2C-(\frac{\alpha}{4h}-2\beta)i},$$ for any $(i,j)\in\Lambda$ with
$i+j\in 2\N+l,l=0,1.$ 
\end{prop}
 
{\bf Proof:} We apply the same argument as in the proof of proposition
\ref{spher-prop}, using lemma \ref{spher-lem-(1,-1)} (which is still true in
characteristic $2$) and lemma \ref{ch2-lem-(0,2)} (in the particular case when
$(V,\tau)$ is the trivial representation of $K$). We will get two limits because
the moves $(i,j)\mapsto (i+1,j-1)$ and $(i,j)\mapsto(i,j+2)$ generate a
sublattice of $\Z^2$ of index $2.$

First, we put $S_\alpha=\{(i,j)\in\Lambda|0\leq i-2j\leq\alpha\}.$
When $0\leq2j\leq i,$ we have
$\big(i,j+2\flf{\frac{i-2j}{4}}\big)\in S_4,$ and by the particular case of
lemma \ref{ch2-lem-(0,2)} when $(V,\tau)$ is the trivial representation of $K,$ we get
$$|c(i,j)-c(i,j+2\flf{\frac{i-2j}{4}})| \leq C'e^{2C-(\frac{\alpha}{4h}-2\beta)i}.$$ 
When $0\leq i\leq 2j$, we have
$\big(i+\clf{\frac{2j-i}{3}},j-\clf{\frac{2j-i}{3}}\big)\in S_3\subset S_4.$ By
lemmas \ref{spher-lem-(1,-1)} we have $$\Big|c(i,j)-c(i+\clf{\frac{2j-i}{3}},j-\clf{\frac{2j-i}{3}})\Big|\leq
C'e^{2C-(\frac{\alpha}{h}-3\beta)\frac{i+j}{3}}.$$ 

Moreover, when $(i,j)\in S_4,$ there exists $k\in\{0,1,2\}$ such that
$$|c(i,j)-c\big(i+k,\frac{1}{2}(i+k)\big)|\leq
C'e^{2C-(\frac{\alpha}{4h}-2\beta)i}.$$ In fact, when $(i,j)\in S_8,$ we first have
$$\tx{max}\Big(|c(i,j)-c(i,j+2)|,|c(i,j)-c(i+1,j-1)|\Big)\leq
C'e^{2C-(\frac{\alpha}{4h}-2\beta)i}.$$ It suffices to show the inequality when
$i-2j\in\{1,2,3,4\},$ by considering the following moves inside $S_8.$ When
$i-2j=1,$ we obtain the inequality by considering
$(2j+1,j)\mapsto(2j+2,j-1)\mapsto(2j+2,j+1).$ When $i-2j=2,$ we consider
$(2j+2,j)\mapsto(2j+4,j-2)\mapsto(2j+4,j+2).$ When $i-2j=3$ or $4,$ use the
moves $(2j+3,j)\mapsto(2j+2,j+1)$ and $(2j+4,j)\mapsto(2j+4,j+2)$ respectively.

In sum, when $i\geq 2j\geq 0, $ there exists $k\in\{0,1,2\},$ such that
\begin{align}|c(i,j)-c(i+k,\frac{1}{2}(i+k))|\leq
C'e^{2C-(\frac{\alpha}{4h}-2\beta)i}\label{zigzag-char2-1};\end{align} when
$2j\geq i\geq j\geq 0,$ there exists $k\in\{0,1,2,3\}$ such that
\begin{align}|c(i,j)-c\big(\flf{\frac{2}{3}(i+j)}+k,\frac{1}{2}(\flf{\frac{2}{3}(i+j)}+k)\big)|\nonumber\\
\leq
C'e^{2C-(\frac{\alpha}{4h}-2\beta)i}\label{zigzag-char2-2}.\end{align} 

Finally the proposition follows from the inequality $$|c(2j,j)-c(2j+4,j+2)|\leq
C'e^{2C-(\frac{\alpha}{4h}-2\beta)2j}.$$
\qed

\begin{prop}\label{nonspher-prop-char2}
Let $\alpha>0$, $\beta\in[0, \frac{\alpha}{8h}),$ and $(V,\tau)$ a non trivial
irreducible unitary representation of $K$. There exists $C'>0$, such that the
following holds. Let $C\in\R_+^*$, $(E,\pi)$ any element in
$\mathcal E_{G,C+\beta\ell}$, and $\xi\in E$, $\eta\in V\otimes E^*$ any
$K$-invariant vectors of norm $1$. We have $$\|c(i,j)\|_V\leq
C'e^{2C-(\frac{\alpha}{4h}-2\beta)i}.$$
\end{prop}

{\bf Proof:} As \eqref{zigzag-char2-1} and \eqref{zigzag-char2-2} in the proof
of the above proposition \ref{spher-prop-char2}, by lemmas \ref{nonspher-lem-(1,-1)} and \ref{ch2-lem-(0,2)}, we have the following
inequalities. When $i\geq 2j\geq 0, $ there exists $k\in\{0,1,2\},$ such that
$$\|c(i,j)-c(i+k,\frac{1}{2}(i+k))\|_V\leq
C'e^{2C-(\frac{\alpha}{4h}-2\beta)i}.$$ When
$2j\geq i\geq j\geq 0,$ there exists $k\in\{0,1,2,3\}$ such that
$$\|c(i,j)-c\big(\flf{\frac{2}{3}(i+j)}+k,\frac{1}{2}(\flf{\frac{2}{3}(i+j)}+k)\big)\|_V\leq
C'e^{2C-(\frac{\alpha}{4h}-2\beta)i}.$$ So it remains to prove that for any
$j\in\N$ we have \begin{align}\|c(2j,j)\|_V\leq
C'e^{2C-(\frac{\alpha}{4h}-2\beta)2j}.\label{zigzag-final}\end{align}

First when $j\in 2\N,$ we know inequality \eqref{zigzag-final} holds.
 In fact, by lemmas
\ref{nonspher-lem-(1,-1)} and \ref{ch2-lem-(0,2)}, when $j\in 2\N$ we have
$$\tx{max}\Big(\|c(2j,j)-c(2j,0)\|_V,\|c(2j,j)-c(3j/2,3j/2)\|_V\Big)\leq  
C'e^{2C-(\frac{\alpha}{4h}-2\beta)2j}. $$ Let $K_1,K_2$ be the subgroups of the
group $K$ as lemma \ref{generated-sg}. By lemmas \ref{generated-sg} and
\ref{invariant} we get inequality \eqref{zigzag-final}.

It remains to show inequality \eqref{zigzag-final} when $j\in 2\N+1.$ We first
have
$$\tx{max}\Big(\|c(2j,j)-c(2j+1,0)\|_V,\|c(2j,j)-c\big(2j-\flf{\frac{j}{2}},j+\flf{\frac{j}{2}}\big)\|_V\Big)$$
$$\leq C'e^{2C-(\frac{\alpha}{4h}-2\beta)2j}.$$ Note that lemma
\ref{nonspher-lem-(1,-1)} is still valid for $i=j-1,$ i.e.
$$\|c\big(2j-\flf{\frac{j}{2}},j+\flf{\frac{j}{2}}\big)-c\big(2j-\flf{\frac{j}{2}}-1,j+\flf{\frac{j}{2}}+1\big)\|_V\leq
C'e^{2C-(\frac{3\alpha}{h}-3\beta)j}.$$ Then we have
$$\|c(2j,j)-c\big(2j-\flf{\frac{j}{2}}-1,j+\flf{\frac{j}{2}}+1\big)\|_V \leq
C'e^{2C-(\frac{\alpha}{2h}-2\beta)2j}.$$

Denote $B_1,B_2$ the image in $K_1$ of
$\begin{pmatrix}\OO^\times&\pi\OO\\\OO&\OO^\times\end{pmatrix}$
and $\begin{pmatrix}\OO^\times&\OO\\\pi\OO&\OO^\times\end{pmatrix}$
respectively, under the group isomorphism $GL_2(\OO)\to K_1.$ We see that
$K_1=(B_1B_2)^2.$ Moreover $c(k'gk'')=\tau(k')c(g)$ for any $k',k''\in K,g\in
G,$ it follows that $c(2j+1,0),
c\big(2j-\flf{\frac{j}{2}},j+\flf{\frac{j}{2}}\big),c\big(2j-\flf{\frac{j}{2}}-1,j+\flf{\frac{j}{2}}+1\big)$
are invariant by $K_2,B_1,B_2$ repectively. By applying lemma \ref{invariant} to
$K=(B_1B_2K_2)^{60},$ we obtain inequality \eqref{zigzag-final} for
$j\in2\N+1.$
\qed 

{\bf Proof of theorem \ref{thm-sp4} when $\tx{char}(F)= 2$:} 
For simplicity we say that an element $g\in G$ is even (resp. odd) when $g\in
KD(i,j)K,i\geq j\geq 0$ and $i+j$ is even (resp. odd). By proposition
\ref{spher-prop-char2}, we see that when $g$ is even
(resp. odd) and tends to infinity, the limit of $e_Ke_ge_K$ exists in
$\CC_{C+\beta\ell}^{\E}(G),$ which we denote by $T_0$ (resp. $T_1$). Let
$\pr=\frac{1}{2}(T_0+T_1).$ 

We have for any $g\in G, $ $e_Ke_g\pr=\pr,$ and thus
$\pr^2=\pr.$ In fact, it suffices to show that 
for any $g\in G$ there exist $\alpha(g),\beta(g)>0$ with $\alpha(g)+\beta(g)=1,$
such that
\begin{align}e_Ke_gT_0=\alpha(g)T_0+\beta(g)T_1,\label{vol-1}\end{align}
and \begin{align} e_Ke_gT_1=\beta(g)T_0+\alpha(g)T_1.\label{vol-2}\end{align}
Let $\alpha(g)$ (resp.
$\beta(g)$) be the volume of the set of elements $(k_1,k_2,k_3,k_4)\in K$ (where $k_i$ are vectors in $F^4$ with norms $\leq 1$)
such that $$\|gk_1\wedge gk_2\|_{\wedge (F^4)}\in q^{2\Z} (\text{resp. } q^{2\Z+1}).$$
We see that for any $k=(k_1,k_2,k_3,k_4)\in K,$ when $i+j\in 2\N$ (resp. $2\N+1$) with $(i,j)\in\Lambda$ and when $\|gk_1\wedge
gk_2\|\geq q^{-2j},$ $gkD(i,j)$ is even exactly when $\|gk_1\wedge gk_2\|\in
q^{2\Z}$ (resp. $q^{2\Z+1}$). Hence we have
$$\lim_{i+j\in2\N(\text{ resp. } 2\N+1),j\to\infty} \text{vol}\{k\in K,
gkD(i,j)\text{ is even}\}=\alpha(g)(\text{ resp. }\beta(g)),$$ and
also $$\lim_{i+j\in2\N(\text{ resp. } 2\N+1),j\to\infty} \text{vol}\{k\in K,
gkD(i,j)\text{ is odd}\}=\beta(g)(\text{ resp. }\alpha(g)).$$ 
And thus equalities \eqref{vol-1} and \eqref{vol-2} follow.

By proposition \ref{nonspher-prop-char2}, for any non trivial irreducible
representation $V$ of $K$ we have $e_K^Ve_gT_0=e_K^Ve_gT_1=0.$ By the same
argument as in the proof of theorem when $\text{char}(F)\neq 2$ in section 2, we
have $$e_g\pr=e_Ke_g\pr=\pr.$$ We complete the proof by taking
$$\pr_n=\frac{1}{2}\big(e_Ke_{D(2\flf{\frac{n}{2}},0)}e_K+
e_Ke_{D(2\flf{\frac{n}{2}}-1,0)}e_K \big),$$ and
$t=\frac{\alpha}{4h}-2\beta.$\qed

\section{Extension to simple algebraic groups of higher split
rank}\label{section-ag}
Let $F$ be a non archimedean local field. This section is dedicated to the proof
the the following theorem, which is theorem \ref{main-thm} in the introduction.
\begin{thm}\label{main}
Let $G$ be a connected almost $F$-simple algebraic group with $F$-split rank $\geq 2$. Then $G(F)$ has
strong Banach property (T).
\end{thm}

We begin the proof with some lemmas. The following lemma is proposition 8.2 in
\cite{bor}. 
\begin{lem}\label{borel}
Let $k$ be a field and $H$ an abelian $k$-group.
Let $\pi:H\to\mathrm{GL}_n$ be a $k$-rational representation. Then $\pi
(H)$ is conjugate over $k$ to some subgroup of the group of diagonal elements in
$\mathrm{GL}_n.$
\end{lem}
The following lemma is a consequence of theorem 7.2 in
\cite{borel-tits}, which is also proposition I.1.6.2 in \cite{mar}.
\begin{lem}\label{borel-tits}
Let $k$ be any field and $G$ a connected almost $k$-simple group with
$k$-split rank $\geq 2.$ Then there exists a $k$-rational group homomorphism
with finite kernel from $SL_3$ or $Sp_4$ to $G.$
\end{lem}
The following lemma is a direct consequence of propositions I.1.3.3 (ii) and
I.1.5.4 (iii), and theorem I.2.3.1 (a) in \cite{mar}.
\begin{lem}\label{mar}
Let $G$ be a simply connected and almost $F$-simple group. Let $S$ be a
maximal $F$-split torus of $G,$ $\Phi(G,S)$ the root system with some ordering
and $\vartheta$ a proper subset of simple roots. Then
there exist two unipotent $F$-subgroups $V_\vartheta,V_\vartheta^-$ of $G,$ and
two $S$-equivariant $F$-isomorphisms $\mathrm{Lie}V_\vartheta\to
V_\vartheta,\mathrm{Lie}V_\vartheta^-\to V_\vartheta^-,$ such that
\begin{itemize}
  \item (i) $\mathrm{Lie}V_\vartheta$ (resp. $\mathrm{Lie}V_\vartheta^-$) is
the direct sum of eigenspaces of positive (resp. negative) roots
which are not integral linear combinations of $\vartheta,$ and
  \item  (ii) $V_\vartheta(F)\cup V_\vartheta^-(F)$ generates $G(F).$
\end{itemize} 
\end{lem}
The next two lemmas
reduce the proof to the simply connected covering of our algebraic group.
\begin{lem}(proposition I.1.4.11 in \cite{mar})\label{covering}
Let $k$ be a field, and let $G$ be connected semisimple $k$-group. Then there
exists a simply connected $k$-group $\tilde G$ and a $k$-isogeny (i.e. surjective
$k$-group homomorphism with finite kernel) from $\tilde G$ to $G.$
\end{lem}
\begin{lem}\label{quotient-finite}
Let $G_1$ be a locally compact group and $G_2$ its quotient by a finite
normal subgroup. Then $G_1$ has strong Banach property (T) if and only if $G_2$
has strong Banach property (T).
\end{lem}
{\bf Proof:} Let $H$ be the kernel of $G_1\to G_2.$ 
 Suppose $G_1$ has strong Banach property (T), and let $\pr_n\in C_c(G_1)$
 be real and self-adjoint elements (otherwise take
 $\pr_n+\bar\pr_n+\pr_n^*+\bar\pr_n^*$) that tends to the idempotent element in
 $\CC_{C+s\ell}^\E(G_1).$ Then $\big(\esp{h\in H}h\big)\pr_n$ tends to a real and self-adjoint (since $H$ is normal) idempotent
element $\pr'$ in $\CC_{C+s\ell}^\E(G_2)$ such that $e_{g}\pr'=\pr'$ for any
$g\in G_2.$ 

On the other direction, if $G_2$ has strong Banach property (T), let
$\pr_n\in C_c(G_2)$ tend to the idempotent element in $\CC_{C+s\ell}^\E(G_2),$
and denote its lifting to $C_c(G_1)$ by $\tilde \pr_n$ (i.e. $\tilde
\pr_n(gh)=\pr_n(g)$ for any $g\in G_1,h\in H$). 
For any
$(E,\pi)\in\E_{G_1,C+s\ell},$ we have
$\pi(\tilde \pr_n)\xi=\pi(\tilde \pr_n)(\esp{h\in H}\pi(h)\xi),$
and thus
$$\|\pi(\tilde\pr_n)-\pi(\tilde\pr_m)\|_{\mathcal
L(E)}\leq\max_{h\in H}\|\pi(h)\|\|\pi(\tilde\pr_n)-\pi(\tilde\pr_m)\|_{\mathcal L(E^H)},$$ where $E^H$
denotes the space of $H$-invariant vectors. We conclude that $\tilde \pr_n$
tends to a real and self-adjoint idempotent element $\pr$ in
$\CC_{C+s\ell}^\E(G_1)$ such that $e_g\pr=\pr$ for any $g\in G_1.$ \qed

{\bf Proof of theorem \ref{main}:} In view of lemmas \ref{covering}
and \ref{quotient-finite}, we can assume $G$ is simply connected (in order to apply lemma \ref{mar} as indicated below).
 By lemma \ref{borel-tits} 
there exist a subgroup $R$ of $G(F)$ and a surjective group
homomorphism $I$ from $SL_3(F)$ or $Sp_4(F)$ to $ R$  with finite
kernel. Let $\rho:F^*\to SL_3(F)$ (resp. $Sp_4(F)$) be the group homomorphism
defined by $$x\mapsto\begin{pmatrix}x&0&0\\0&1&0\\0&0&x^{-1}\end{pmatrix}\text{(resp.}
\begin{pmatrix}x&&&\\&1&\\&&1&\\&&&x^{-1}\end{pmatrix})$$ for any $x\in F$, and
let $a=I\circ\rho(\pi)$, where $\pi$ is a uniformizer of $F$.
By lemma \ref{borel}, the set of eigenvalues of $Ad(a)$ is a subset of
$\pi^{\Z}$ which contains $\{1\}$ as a proper subset. Let $S$ be a maximal
$F$-split torus of $G$ whose $F$ points contains $a$. We can choose an
ordering of $\Phi(S,G)$ such that $|\chi(a)|\leq1$ for any simple root $\chi$. Let $\vartheta$ be the proper
subset of simple roots $\chi$ such that $|\chi(a)|=1$, and let $V_\vartheta,
V^-_\vartheta$ be as in lemma \ref{mar}. 

For simplicity denote $G(F)$ and $V_\vartheta(F),
V^-_\vartheta(F)$ by $G$ and $V_\vartheta,
V^-_\vartheta$ from now on. Let $\|\cdot\|$ be
the norm on $\lie G$ defined w.r.t.
some $F$-basis.
 Let $\ell'$ be the
length function on $G$ defined by $$\ell'(g)=\log \|Ad(g)\|_{\mathrm{End}(\lie
G)}.$$ 

Note that for any length function $\tilde \ell$ on $SL_3(F)$ or $Sp_4(F),$ there exist $\kappa\in\R_+^*$ such that $\tilde\ell\leq\kappa\ell,$ where $\ell$ is the length function on $SL_3(F)$ or $Sp_4(F)$ defined in section \ref{section-sp4}. In fact for $SL_3(F),$ let $K$ be the compact generating set $$\{SL_3(\OO),SL_3(\OO)\pi\begin{pmatrix}
\pi^{-3}&&\\&1&\\&&1\end{pmatrix}SL_3(\OO),SL_3(\OO)\pi^2\begin{pmatrix}
\pi^{-3}&&\\&\pi^{-3}&\\&&1\end{pmatrix}SL_3(\OO)\}.$$ 
Then we have $\ell(g)\geq\min\{n:g\in K^n\}$ (note that it holds for diagonal elements and $K$ is a $SL_3(\OO)$ bi-invariant set). Therefore, we have
$$\Big(\max_{g\in K}{\tilde\ell(g)}\Big)\ell(g)\geq\tilde \ell(g).$$
It can be shown for $Sp_4(F)$ using the same argument by replacing the compact generating set $K$ by  
$$\{Sp_4(\OO),Sp_4(\OO)\begin{pmatrix}
\pi^{-1}&&&\\&1&&\\&&1&\\&&&\pi\end{pmatrix}Sp_4(\OO),Sp_4(\OO)\begin{pmatrix}
\pi^{-1}&&&\\&\pi^{-1}&&\\&&\pi&\\&&&\pi\end{pmatrix}Sp_4(\OO)\}.$$ 

Let $\E$ be a class of Banach spaces of type $>1$ stable under duality
and complex conjugation. Let $s,t,C,C'\in\R_{+}^{*}, \pro\in \CC^{\E}_{s\kappa\ell+C}(R),\pro_{m}\in C_{c}(R)$ verify the conditions (i) and
(ii) of theorem \ref{thm-sl3} if $R$ is isogenous to $SL_3(F)$, or of theorem
\ref{thm-sp4} if $R$ is isogenous to $Sp_4(F)$, where $\kappa\in\R_{+}^{*}$
such that $\ell'|_{R}\leq\kappa\ell$ (in view of lemma \ref{quotient-finite}). 
Let $U$ be an open compact subgroup of $G$ and $f=\frac{e_U}{vol(e_U)}$. Then
to establish that $G$ has strong Banach property (T) it suffices to
show that if $s$ is small enough the series $\pro_{m}f\in C_{c}(G)$ converges
in $\CC^{\E}_{s\ell'+C}(G)$ to a self adjoint idempotent $\pro'$ such that
for any $(E,\pi)\in\E_{G,s\ell'+C}$, the image of
$\pi(\pro')$ consists of all $G$-invariant vectors of $E.$
First it is clear that the series $\pro_{m}f$ is a Cauchy series in $\CC^{\E}_{s\ell'+C}(G)$ and we note $\pro'$
its limit (in fact $\pro$ is a multiplier of $\CC^{\E}_{s\ell'+C}(G)$ and
$\pro'=\pro f$). Let $(E,\pi)\in \E_{G,s\ell'+C}$. It is obvious that
$\pi(\pro')$ acts by identity over any $G$-invariant vector. It remains to show that for any $x\in E$, $\pi(\pro')x$
is $G$-invariant (in fact it follows that $\pro'=f^{*}\pro'=f^{*}\pro f$, so
$\pr'$ is self-adjoint).
In view of statement (ii) of lemma \ref{mar}, it suffices to show that $\pi(\pro')x$ is $V_\vartheta$-invariant and
$V^-_\vartheta$-invariant. 

We first show that $\pi(\pro')x$ is $V_\vartheta$-invariant. Let $E:\lie V_\vartheta\to V_\vartheta$
be as in lemma \ref{mar}.
We know that $\pi(\pro')x$ is fixed by $R$, then in particular by $a$. It
suffices to show that for any $Y\in\lie V_\vartheta,$
$$\pi(E(Y))\pi(\pro')x-\pi(\pro')x=\pi(E(Y))\pi(a^{-n})\pi(\pro')x-\pi(a^{-n})\pi(\pro')x$$
$$=\pi(a^{-n})\big(\pi(a^{n}E(Y)a^{-n})-1\big)\pi(\pro')x=\pi(a^{-n})\big(\pi(E(Ad(a^{n})Y))-1\big)\pi(\pro')x$$
tends to $0$ when $n\in\N$ tends to infinity.

Let $Y=\sum_{\lambda\in\Lambda} Y_\lambda$ be the decomposition of $Y$ under the adjoint action of $a$ in $\lie
V_\vartheta$, where
$\Lambda\subset F$ denotes the set of eigenvalues
of the action.
Due to the way $\vartheta$ is chosen, the eigenvalues of $Ad(a)|_{\lie
V_\vartheta}$ are all of the form $\pi^{\N^*}.$ Since
$U$ is an open subgroup of $G$, there exists $r>0$ such that when $Y'\in V$ and
$\|Y'\|\leq r,$ we have $E(Y')\in U$. For any $n\in\N,$ we put
$$m=\flf{n\kappa^{-1}\log\min_{\lambda\in\Lambda}|\lambda|^{-1}+\kappa^{-1}\log(r/\max_{\lambda\in\Lambda}
\|Y_\lambda\|)}.$$ When $n$ is big enough such that $m>0$, we have
$$\big(\pi(E(Ad(a^{n})Y))-1\big)\pi(\pro_{m}f)x$$ $$=\int_{R}
\pro_{m}(g)\pi(g)\big(\pi(E(Ad(g^{-1}a^{n})Y))-1\big)\pi(f)x dg.$$
 When $\ell(g)\leq m,$ 
 we have $$\|Ad(g^{-1}a^{n})Y\|\leq
 e^{\ell'(g)}\max_{\lambda\in\Lambda}|\lambda|^{n}\|Y_\lambda\|_{\lie V_\vartheta}\leq r,$$ and hence
 $$\big(\pi(E(Ad(g^{-1}a^{n})Y))-1\big)\pi(f)x=0.$$ Therefore we have
 $$\pi(a^{-n})\big(\pi(E(Ad(a^{n})Y))-1\big)\pi(\pro_{m}f)x=0$$ when $n$ is
 big enough.          

On the other hand for any $n\in\N,$
we always have $$Ad(a^{n})Y=\sum_{\lambda\in\Lambda}
\lambda^{n} Y_\lambda\in \bigoplus_{\lambda\in\Lambda} \OO Y_\lambda.$$ Hence
$$\|\pi(a^{-n})\big(\pi(E(Ad(a^{n})Y))-1\big)\pi(\pro'-\pro_{m}f)x\|_{E}$$
$$\leq e^{C+s\ell'(a)n} (1+C'')\|\pi(\pro'-\pro_{m}f)x\|_{E},$$ where
$$C''=\sup_{t_\lambda\in \OO}\|\pi(E(\max_{\lambda\in\Lambda} t_\lambda
Y_\lambda))\|_{\mathcal{L}(E)}<\infty$$ depends only on $Y$. But 
$$\|\pi(\pro'-\pro_{m}f)x\|_{E}\leq C' e^{2C-tm} \|\pi(f)x\|_{E} $$ by
statement (ii) of theorem \ref{thm-sl3} if $R$ is isogenous to $SL_3(F)$, or of
theorem \ref{thm-sp4} if $R$ is isogenous to $Sp_4(F)$ (we recall that $C'$ and
$t$ are the constants of theorem \ref{thm-sl3} and theorem \ref{thm-sp4}). In
total, when $n$ is big enough
$$\|\pi(a^{-n})\big(\pi(E(Ad(a^{n})Y))-1\big)\pi((\pro'-\pro_{m}f))x\|_{E}$$
$$\leq e^{C+s\ell'(a)n} (1+C'') C' e^{2C-tm} \|\pi(f)x\|_{E},$$ and if
$$s<\frac{t}{\kappa\ell'(a)}\log\min_{\lambda\in\Lambda}|\lambda|^{-1},$$
it tends to $0$ when $n$ tends to infinity. 

 We prove
$\pi(\pro')x$ is $V_\vartheta^-$-invariant by exactly the same argument (with
$a$ replaced by $a^{-1}$ and the ordering of $\Phi(S,G)$ by its inverse, i.e. the ordering such that $|\chi(a^{-1})|\leq 1$ for any simple root $\chi$).
\qed

\end{document}